\documentclass[11pt,leqno]{article} 

\usepackage{fullpage}

\usepackage{amssymb} 
\usepackage{amsmath} 
\date{}
\title{Algebraic vs physical N=6 3-algebras} 
\author{{\sc Nicoletta Cantarini}\thanks{Dipartimento di Matematica, Universit\`a di Bologna, Bologna, Italy -
Supported in part by
Progetto FIRB ``Perspectives in Lie theory"}
\and\setcounter{footnote}{6}
{\sc Victor G.\ Kac}
\thanks{Department of Mathematics, MIT, Cambridge, Massachusetts 02139, USA}}

\newtheorem{theorem}{Theorem}[section] 
\newtheorem{lemma}[theorem]{Lemma} 
 
\newtheorem{proposition}[theorem]{Proposition} 
\newtheorem{definition}[theorem]{Definition} 
\newtheorem{remark}[theorem]{Remark}

\newtheorem{example}[theorem]{Example}

\def\Z{\mathbb{Z}}
\def\R{\mathbb{R}}
\def\ZZ{\mathbb{Z}}
\def\g{\mathfrak{g}}

\def\C{\mathbb{C}}
\def\F{\mathbb{F}}

\def\ZZ{\mathbb{Z}}

\def\0{\bar{0}}
\def\1{\bar{1}}

\numberwithin{equation}{section}

\makeatletter

\def\enumerate{%
  \ifnum \@enumdepth >\thr@@\@toodeep\else
    \advance\@enumdepth\@ne
    \edef\@enumctr{enum\romannumeral\the\@enumdepth}%
      \list
        {\csname label\@enumctr\endcsname}%
        {\usecounter\@enumctr
          \addtolength{\leftmargin}{-\leftmargin}
          \settowidth{\labelwidth}{(99)}
          \itemindent = \labelwidth
          \addtolength{\itemindent}{\labelsep}
        \listparindent=1em      
          \def\makelabel##1{{##1}\hfill}
          }%
  \fi}

\makeatother

\begin{document} 
\maketitle 
\begin{abstract} In our previous paper we classified linearly compact algebraic simple $N=6$ $3$-algebras. In the present
paper we classify their ``physical" counterparts, which actually appear
in the $N=6$ supersymmetric 3-dimensional Chern-Simons theories.
\end{abstract} 

\setcounter{section}{-1} \label{sec:intro}
\section{Introduction}
In a series of papers on $N$-supersymmetric $3$-dimensional Chern-Simons gauge theories various types of $3$-algebras have naturally
appeared (\cite{G}, \cite{BL1}, \cite{BL2}, \cite{ABJM}, \cite{MFM}, \cite{ST}, \cite{BB}, \dots).

Recall that a $3$-algebra is a vector space $A$ over $\C$ with a $3$-bracket $A^{\otimes 3}\rightarrow A$, $a\otimes b\otimes c\mapsto
[a,b,c]$. If this bracket is linear in all arguments, the $3$-algebra $A$ is called {\it algebraic}. If this $3$-bracket is linear in the first and the third argument, but anti-linear in the second argument (i.e.\ $[\lambda a, b, c]=[a,b,\lambda c]=\lambda [a,b,c]$, but $[a, \lambda b,c]=
\bar{\lambda}[a,b,c]$, where $\bar{\lambda}$ is the complex conjugate of $\lambda\in\C$), then $A$ is called a {\it physical} $3$-algebra. 

If $C$ is an anti-linear involution of the vector space $A$  (i.e. $C^2=1$ and $C(\lambda a)=\bar{\lambda}C(a)$), then an algebraic $3$-algebra
 $A$ with the $3$-bracket $[a,b,c]$ is converted to a physical $3$-algebra with 
the $3$-bracket
\begin{equation}
[a,b,c]_{ph, C}=[a, C(b),c],
\label{0.1}
\end{equation}
which we denote by $A_{ph, C}$, and viceversa.

An algebraic or physical $3$-algebra $A$ is called $N=6$ $3$-algebra if it satisfies the following two axioms, cf.\ \cite{BL1} and \cite{CantaK6} ($a,b,c,x,y,z\in A$):
$$
[a,b,c]=-[c,b,a]~~~~({\mbox{anti-commutativity}})$$
$$[a,b,[x,y,z]]=[[a,b,x],y,z]-[x,[b,a,y],z]+[x,y,[a,b,z]]~~ ({\mbox{fundamental identity}}).$$
In our previous paper \cite{CantaK6} we classified all linearly compact simple $N=6$ algebraic $3$-algebras. It is easy to see that if
$C$ is an anti-linear involution of the $3$-bracket $[a,b,c]$ of an algebraic $N=6$ $3$-algebra, then $[a,b,c]_{ph,C}$, defined by
(\ref{0.1}), is a $3$-bracket of a physical $N=6$ $3$-algebra. However not all physical $N=6$ 3-algebras are obtained from the algebraic ones in 
this way, and we call them {\it trivially related} if they are.
Indeed, it turned out that the theory of physical $N=6$ $3$-algebras is richer than that of algebraic ones. On the other hand, the method 
used in the classification of the former is similar to that of the latter.

Recall that the space $M_{m,n}(\C)$ of all $m\times n$ matrices over $\C$ carries a structure of an algebraic $N=6$ $3$-algebra, given by
\begin{equation}
[a,b,c]=ab^tc-cb^ta,
\label{0.2}
\end{equation}
where $a\mapsto a^t$ is the transpose, and this $3$-algebra is simple. It is denoted by $A^3(m,n;t)$ in \cite{CantaK6}.

On the other hand, the space $M_{m,n}(\C)$ carries the following structure of physical $N=6$ $3$-algebra, given by
\begin{equation}
[a,b,c]=aC_{p,q}(b)^tc-cC_{p,q}(b)^ta,
\label{0.3}
\end{equation}
where $C_{p,q}$ is the anti-linear involution of the space $M_{m,n}(\C)$, defined by
$$C_{p,q}(a)=S^m_q\bar{a}S^n_p,
\, {\mbox{where}}\, S^m_q=diag(I_q, -I_{m-q}), ~0\leq q\leq m, ~0\leq p\leq n$$
(and $I_q$ stands for the $q\times q$ identity matrix).
This is a simple physical $N=6$ $3$-algebra, denoted by $A^3(m,n;t)_{ph, C_{p,q}}$, which is, by definition, trivially related to
the algebraic $N=6$ $3$-algebra $A^3(m,n;t)$; for $p=m$, $q=n$ it was introduced in \cite{BL2}.

The simplest example of a simple physical $N=6$ $3$-algebra, which is not trivially related to any algebraic one, is obtained by endowing the
space $M_{n,n}(\C)$ with the 3-bracket
\begin{equation}
[a,b,c]=\pm i(a\bar{b}c-c\bar{b}a)
\label{0.4}
\end{equation}
(different signs give non-isomorphic 3-algebras). It is denoted by $A^3(n)_{\pm}$.

An interesting phenomen here is that physical $N=6$ 3-algebras, trivially related to non-isomorphic algebraic 
$N=6$ 3-algebras, may be isomorphic (see Remark \ref{phenomenon}).

The first main result of the paper is the classification of finite-dimensional simple physical
$N=6$ 3-algebras, given by the first part of Theorem \ref{main}. The answer consists of four series: the two series of physical 3-algebras 
$A^3(m,n;t)_{ph, C_{p,q}}$ and $A^3(n)_{\pm}$ of type $A$, mentioned above, and two more series of type $C$. One of the series
of type $C$, denoted by $C^3(2n)_{ph, \pm C_{p}}$, where
$$C_{n-p}(u)=\bar{u}\,diag(S^n_p, S^n_p)~~(0\leq p\leq n),$$
introduced for $p=n$ in \cite{BL2}, is trivially related to the
algebraic 3-algebra $C^3(2n)$, introduced in \cite{CantaK6}, and the second series of type $C$, denoted by $C^3(2n,iS_n^{2n}, \pm i)$,
introduced in Example \ref{C(n)} of the paper, is not trivially related to any algebraic $N=6$ $3$-algebra.

The second result of the paper is the classification of infinite-dimensional linearly compact simple physical $N=6$ 3-algebras, given by 
the second part of Theorem \ref{main}. The third part of this theorem describes which of the examples are not trivially related
to any algebraic $N=6$ 3-algebra. The most interesting of them is the family of physical $N=6$ 3-algebras $W^3_{\beta}(\varphi)_{\pm}$, 
depending on a complex parameter $\beta$ of modulus 1, $\beta\neq\pm 1$, which is not trivially related to any algebraic $N=6$
3-algebra.

As in \cite{CantaK6}, in order to achieve the classification, we establish in Theorem \ref{th:2.3} a bijection between the isomorphism
classes of physical $N=6$ 3-algebras $A$ with zero center and the isomorphism classes of pairs $(\g, \sigma)$, where 
$\g=\oplus_{j\in\Z}\g_j$ is a $\Z$-graded Lie superalgebra with a consistent short grading, satisfying certain
additional properties (which in case of simple $A$ mean that $\g$ is simple) and $\sigma$ is an anti-linear graded conjugation of $\g$.
The only difference with the algebraic case, considered in \cite{CantaK6}, is that $\sigma$ in \cite{CantaK6} is a linear graded
conjugation (graded conjugation means that $\sigma(\g_j)=\g_{-j}$ and $\sigma^2|_{\g_k}=(-1)^k$).

Thus, the classification in question reduces to the classification of anti-linear graded conjugations of simple linearly compact
Lie superalgebras with consistent short $\Z$-gradings (such $\Z$-graded Lie superalgebras have been classified in \cite{CantaK6}).
This is technically the most difficult part of the paper.

Note a relation of $N=6$ 3-algebras to Jordan triple systems (called also Jordan 3-algebras), introduced by Jacobson \cite{J},
as subspaces of an associative algebra, closed under the 3-bracket $[a,b,c]=abc+cba$, and axiomatized in \cite{Mey}. Doing the same in an associative superalgebra
we arrive at the following:
\begin{definition}\label{N6} A Jordan $3$-superalgebra is a vector superspace endowed with an algebraic 
$3$-bracket $[\cdot,\cdot,\cdot]$ satisfying
the following axioms:
\begin{enumerate}
\item[$(a)$] $[a,b,c]=(-1)^{p(a)p(b)+p(a)p(c)+p(b)p(c)}[c,b,a]$
        \item[$(b)$] 
        $[a,b,[x,y,z]]=[[a,b,x],y,z]-(-1)^{p(x)(p(a)+p(b))+p(a)p(b)+p(a)}[x,[b,a,y],z]+$
        \end{enumerate}
\hfill $(-1)^{(p(x)+p(y))(p(a)+p(b))}[x,y,[a,b,z]]$
\end{definition}
Obviously, a purely even Jordan $3$-superalgebra is a Jordan triple system, and a purely odd Jordan 3-superalgebra becomes
an $N=6$ 3-algebra after reversing the parity.
We are planning to classify simple linearly
compact Jordan 3-superalgebras in a subsequent paper \cite{CantaK7}.

We are grateful to the referee of our paper  \cite{CantaK6} who suggested us to study physical $N=6$ 3-algebras. 
\section{Examples of physical $\mathbf{N=6}$ 3-algebras}
%
%
%
First, as has been mentioned in the Introduction, we have the following:
\begin{remark}\label{fromalgtophys}\em Let $(A,[\cdot,\cdot,\cdot])$ be an algebraic $N=6$ $3$-algebra over $\C$ and let $C$ be 
an anti-linear involution of the 3-algebra $A$. Then it is straightforward to check that $A$ with the 3-bracket $[a,b,c]_{ph,C}=[a,C(b),c]$ is a physical $N=6$ $3$-algebra.
We shall denote it by $A_{ph,C}$.
\end{remark}
Note also the following:
\begin{remark}\label{scalars}\em
If $(A,[\cdot, \cdot, \cdot])$ is an algebraic $N=6$ 3-algebra, then for every $\lambda\in\C$, $\lambda\neq 0$,
$(A,[\cdot, \cdot, \cdot]_\lambda)$ where $[a,b,c]_\lambda=\lambda[a,b,c]$, is an algebraic $N=6$ 
3-algebra isomorphic to $(A,[\cdot, \cdot, \cdot])$,
the isomorphism being the map $f_\alpha: A \rightarrow A$, $a\mapsto\alpha a$, for $\alpha\in\C$ such that $\alpha^2=\lambda^{-1}$.

Differently from the algebraic case,  if $(A,[\cdot, \cdot, \cdot])$ is a physical  $N=6$ 3-algebra and $\lambda\in\C$, $\lambda\neq 0$,
$(A,[\cdot, \cdot, \cdot]_\lambda)$ is a physical $N=6$ 3-algebra if and only if $\lambda\in\R$. Besides, the map $f_{\alpha}$ is
an isomorphism of physical $N=6$ 3-algebras if and only if $\lambda\in\R^{>0}$.
\end{remark}

Now we introduce all the examples that appear in our classifications of algebraic and physical $N=6$ 3-algebras. The fact that they indeed satisfy
the fundamental identity either can be checked directly or follows automatically from Theorem \ref{tel}.

\begin{example}\label{sl(2,2)}\em
We denote by $A^3(m,n;t)$ the algebra $M_{m,n}(\C)$ with $3$-bracket:
$$[a,b,c]=ab^tc-cb^ta.$$
This is an algebraic $N=6$ 3-algebra \cite[Example 1.2]{CantaK6}.
The complex conjugation $C_0$ of matrices, 
$$C_0: M_{m,n}(\C)\rightarrow M_{m,n}(\C)$$
$$(a_{ij})\mapsto (\bar{a}_{ij}),$$
is an anti-linear involution of the $3$-algebra $A^3(m,n;t)$, hence, by Remark \ref{fromalgtophys}, 
$M_{m,n}(\C)$ with $3$-bracket:
$$[a,b,c]=a\bar{b}^tc-c\bar{b}^ta$$
is a physical $N=6$ 3-algebra that we shall denote by $A^3(m,n;t)_{ph,C_0}$.

More generally, let  ${}^*: M_{m,n}(\C) \rightarrow M_{n,m}(\C)$ be an anti-linear involutive map satisfying the 
following property:
\begin{equation}\label{*condition}
(ab^*c)^*=c^*ba^*,~~a,b,c\in M_{m,n}(\C).
\end{equation} 
Then $M_{m,n}(\C)$ with the $3$-bracket 
\begin{equation}
[a,b,c]=a{b}^*c-c{b}^*a
\label{antitranspose}
\end{equation}
 is a 
physical $N=6$ 3-algebra.

For example,
 let
$$S_p^n=diag(I_p,-I_{n-p})\in GL_n(\C),\,  ~0\leq p\leq n,$$ and define the map 
\begin{equation}\label{signature}
\varphi_{p,q}: M_{m,n}(\C) \rightarrow M_{n,m}(\C), ~~u~\longmapsto~ S_p^n\bar{u}^tS_q^m.
\end{equation} 
Then $\varphi_{p,q}$ satisfies property (\ref{*condition}), hence  $M_{m,n}(\C)$ with the corresponding 3-bracket
(\ref{antitranspose}) is a physical $N=6$ 3-algebra. Again, we can obtain this 3-algebra 
from the algebraic $N=6$ 3-algebra $A^3(m,n;t)$, as explained in Remark \ref{fromalgtophys}, using the anti-linear involution
$$C_{p,q}: M_{m,n}(\C) \rightarrow  M_{m,n}(\C),~~
u\mapsto S_q^m\bar{u}S_p^n.$$
We shall therefore denote the physical $N=6$ 3-algebra 
$M_{m,n}(\C)$ with the $3$-bracket:
$[a,b,c]=a\varphi_{p,q}(b)c-c\varphi_{p,q}(b)a$ by $A^3(m,n;t)_{ph,C_{p,q}}$.
Note that $C_0=C_{0,0}=C_{n,m}$ and $-C_0=C_{n,0}=C_{0,m}$.

More generally, if $A\in GL_n(\C)$ and $B\in GL_m(\C)$ are such that $A=\lambda \bar{A}^t$, $B=\lambda \bar{B}^t$, for some $\lambda\in\C$, 
then the map $\varphi_{A,B}: M_{m,n}(\C) \rightarrow M_{n,m}(\C)$ 
defined by: $\varphi_{A,B}(b)=A\bar{b}^tB^{-1}$, satisfies property (\ref{*condition}), hence $M_{m,n}(\C)$ with the corresponding 3-bracket
(\ref{antitranspose}) is a physical $N=6$ 3-algebra. Note that $\varphi_{p,q}=\varphi_{S_p^n, S_q^m}$.

\end{example}
\begin{lemma}\label{sl(m,n)}
Let ${}^*: M_{m,n}(\C) \rightarrow M_{n,m}(\C)$ 
be defined by: $b^{*}=A\bar{b}^tB^{-1}$ for some matrices
$A\in GL_n(\C)$ and $B\in GL_m(\C)$ such that $A=\lambda \bar{A}^t$, $B=\lambda \bar{B}^t$, with $\lambda\in\C$, $|\lambda|=1$. Then the 
physical $N=6$ 3-algebra with the corresponding bracket 
(\ref{antitranspose}) is
isomorphic to $A^3(m,n;t)_{ph,C_{p,q}}$,
$p$ and $q$ being the numbers of positive eigenvalues of the hermitian matrices
$\lambda^{-1/2}A$ and $\lambda^{-1/2}B$, respectively.
\end{lemma}
{\bf Proof.} It is immediate to check that, under the hypotheses of the lemma, the matrices 
$A'=\lambda^{-1/2}A$ and $B'=\lambda^{-1/2}B$ are hermitian. It follows that there exist unique matrices $h\in GL_n(\C)$
and $k\in GL_m(\C)$ such that $A'=hS_p^n\bar{h}^t$ and $B'=kS_q^m\bar{k}^t$, where $p$ (resp.\ $q$) is the number of positive eigenvelues
of $A'$ (resp.\ $B'$).
Now consider the map $f: M_{m,n}(\C) \rightarrow M_{m,n}(\C)$
defined by $f(u)=kuh^{-1}$. For $a,b,c\in M_{m,n}(\C)$, let $[a,b,c]_{p,q}=a\varphi_{p,q}(b)c-c\varphi_{p,q}(b)a$
where $\varphi_{p,q}$ is map (\ref{signature}), and $[a,b,c]_*=ab^*c-cb^*a$. Then we have:
$$f([a,b,c]_{p,q})=f(a\varphi_{p,q}(b)c-c\varphi_{p,q}(b)a)=
ka\varphi_{p,q}(b)ch^{-1}-kc\varphi_{p,q}(b)ah^{-1}
=(kah^{-1})(hS_p^n\bar{b}^tS_q^mk^{-1})$$ $$(kch^{-1})
-(kch^{-1})(hS_p^n\bar{b}^tS_q^mk^{-1})(kah^{-1})
=f(a)(A(\bar{h}^t)^{-1}\bar{b}^t\bar{k}^tB^{-1})f(c)-f(c)(A(\bar{h}^t)^{-1}\bar{b}^t\bar{k}^tB^{-1})f(a)$$
$$=f(a)f(b)^*f(c)-f(c)f(b)^*f(a)=[f(a), f(b), f(c)]_*,~~~~~~~~~~~~~~~~~~~~~~~~~~~~~~~~~~~~~~~~~~~~~~~~~~~~~~~~~$$
and this shows that the $3$-brackets $[\,,\,,]_{p,q}$ and $[\,,\,,]_*$ are isomorphic.
\hfill$\Box$

\begin{remark}\label{phenomenon}\em
Recall \cite{CantaK6} that for $m=2h$ and $n=2k$, the following algebraic $N=6$ 3-bracket is defined on $M_{m,n}(\C)$:
$$[a,b,c]=ab^{st}c-cb^{st}a,$$
where the map $st: M_{m,n}(\C) \rightarrow M_{n,m}(\C)$,
is defined by:
 $a\mapsto a^{st}:=J_{2k}a^tJ_{2h}^{-1}$, with 
 $J_{2k}=
 \left(\begin{array}{cc}
 0 & I_k\\
 -I_k & 0
 \end{array}
 \right)$.
 In \cite{CantaK6} we denoted $M_{m,n}(\C)$ with this 3-bracket by $A^3(m,n;st)$. 
 By Remark \ref{fromalgtophys}, we can associate to $A^3(m,n;st)$ the physical $N=6$ 3-algebra $A^3(m,n;st)_{ph,C_0}$ by defining on $M_{m,n}(\C)$ the following 3-bracket:
 $$[a,b,c]=a\bar{b}^{st}c-c\bar{b}^{st}a.$$
 Since $\bar{J}_{2k}^t=-J_{2k}$, by Lemma \ref{sl(m,n)},   $A^3(2h,2k;st)_{ph,C_0}$
 is isomorphic to $A^3(2h,2k;t)_{ph,C_{k,h}}$.
\end{remark}

\begin{example}\label{sl(n,n)}\em 
Let ${\cal A}$ be an associative algebra over $\C$ and let ${}^\dag: {\cal A}\rightarrow {\cal A}$ be an anti-linear map such that:
\begin{equation}
\label{*antilinear}
(ab^\dag c)^\dag=a^\dag bc^\dag, ~~a,b,c\in {\cal A}.
\end{equation}
Then ${\cal A}$ with the 3-bracket
\begin{equation}
\label{newforsl(n,n)}
[a,b,c]=i(ab^\dag c-cb^\dag a)
\end{equation}
is an $N=6$ physical 3-algebra.
For example, if ${\cal A}=M_{n,n}(\C)$ and $\dag$ is the complex conjugation of matrices, then the corresponding bracket (\ref{newforsl(n,n)}) defines on ${\cal A}$ a physical $N=6$
3-algebra structure that we shall denote by $A^3(n)_{+}$. Similarly, we shall denote by $A^3(n)_-$ the physical $N=6$
3-algebra ${\cal A}$ with the negative of the 3-bracket (\ref{newforsl(n,n)}).

More generally, for $A\in GL_n(\C)$ consider the following anti-linear map 
$$\psi_A: M_{n,n}(\C)\rightarrow M_{n,n}(\C),~ ~\psi_A(u)=A\bar{u}\bar{A}.$$ Then $\psi_A$ satisfies property (\ref{*antilinear}), hence $M_{n,n}(\C)$ with the corresponding
bracket (\ref{newforsl(n,n)})
is a physical $N=6$ 3-algebra.
\end{example}

\begin{lemma}\label{lemmasl(n,n)} For every $A\in GL_n(\C)$, $M_{n,n}(\C)$ with the 3-bracket (\ref{newforsl(n,n)}) associated to the map
$\psi_A$ is isomorphic to $A^3(n)_+$.
\end{lemma}
{\bf Proof.} Let us denote by $[\cdot,\cdot,\cdot]_A$ bracket (\ref{newforsl(n,n)}) corresponding to the map $\psi_A$, $A\in GL_n(\C)$,
i.e., $[a,b,c]_A=i(a\psi_A(b)c-c\psi_A(b)a)$
and let $[\cdot,\cdot,\cdot]$ be the 3-bracket in $A^3(n)_+$, i.e., $[a,b,c]=i(a\bar{b}c-c\bar{b}a)$. 
Now let us write $A=hk$ for some $h,k\in GL_n(\C)$ and consider the map $f: M_{n,n}(\C) \rightarrow M_{n,n}(\C)$ defined by
$f(u)=\bar{k}uh$. Then we have:
$f([a,b,c]_A)=i\bar{k}(a\psi_A(b)c-c\psi_A(b)a)h=
i\bar{k}(aA\bar{b}\bar{A}c-cA\bar{b}\bar{A}a)h=i((\bar{k}ah)h^{-1}A\bar{b}\bar{A}(\bar{k})^{-1}(\bar{k}ch)-
(\bar{k}ch)h^{-1}A\bar{b}\bar{A}(\bar{k})^{-1}(\bar{k}ah))=
i(f(a)k\bar{b}\bar{h}f(c)-f(c)k\bar{b}\bar{h}f(a))=i(f(a)\overline{f(b)}f(c)-f(c)\overline{f(b)}f(a))=[f(a), f(b), f(c)]$,
which proves that the 3-brackets $[\cdot,\cdot,\cdot]_A$ and $[\cdot,\cdot,\cdot]$ are isomorphic.
\hfill$\Box$

\medskip

\begin{example}\label{C(n)}\em
Let us consider the map $\psi: M_{1,2n}(\C)\rightarrow M_{2n,1}(\C)$, defined by:
$\psi(X ~Y)=(Y~ -X)^t$, for $X,Y\in M_{1,n}$. Then $M_{1,2n}$ with the $3$-bracket
\begin{equation}
[a,b,c]=-a\bar{b}^tc+c\bar{b}^ta-c\psi(a)(\psi(\bar{b}))^t
\label{C(n)product}
\end{equation}
is a physical $N=6$ $3$-algebra, which we denote by $C^3(2n)_{ph, C_0}$, since it can be obtained from \cite[Example 1.4]{CantaK6}
using Remark \ref{fromalgtophys} with $C_0$ equal to the complex conjugation. Indeed $C^3(2n)$ denotes the algebra $M_{1,2n}(\C)$
with 3-bracket 
$$
[a,b,c]=-a{b}^tc+c{b}^ta-c\psi(a)(\psi({b}))^t
$$
and the complex conjugation of matrices is an anti-linear involution of $C^3(2n)$ with respect
to this bracket, since $\overline{\psi(x)}=\psi({\bar{x}})$.

More generally, let $H\in M_{2n,2n}(\C)$ be a symplectic matrix for the bilinear form with the matrix $J_{2n}$. Then
$M_{1,2n}$ with $3$-bracket
\begin{equation}
[a,b,c]_{H,\alpha}=-\alpha aH\overline{b}^tc+\alpha cH\overline{b}^ta-\alpha^{-1}c\psi(a)(\psi(\overline{b}H^t))^t,
\label{generalizedC(n)product}
\end{equation} 
where either $\alpha\in\R$ and $H$ is hermitian, or $\alpha\in i\R$ and $H$ is anti-hermitian (i.e., $H=-\bar{H}^t$), 
is a physical $N=6$ $3$-algebra  which we denote by $C^3(2n, H; \alpha)$. 
\end{example}

Let $S^n_p\in GL_n(\C)$ be the matrix defined in Example \ref{sl(2,2)}. We set 
$H^{2n}_{p}=diag(S^n_p, S^n_p)$.

\begin{lemma}\label{fromFassbenderIkramov}
Let $H\in M_{2n,2n}(\C)$ be a symplectic matrix (for the bilinear form with the matrix $J_{2n}$). Then
\begin{itemize}
\item[(a)] If $H$ is hermitian, then there exists a symplectic matrix $V$ such that
$H=VH^{2n}_p\overline{V}^t$
where $2p$ is the number of positive eigenvalues of $H$.
\item[(b)] If $H$ is anti-hermitian, then there exists a symplectic matrix $V$ such that
$H=iVS^{2n}_n\overline{V}^t$.
\end{itemize}
\end{lemma}
{\bf Proof.} $(a)$ follows immediately from \cite[Proposition 3]{FI}. In order to prove $(b)$ one can use the same
argument as in \cite[Proposition 3]{FI}. Namely, suppose that $H$ is symplectic and anti-hermitian and set $K=iH$. Then
$K$ is hermitian and satisfies relation $K^tJ_{2n}K=-J$. Let $v$ be an eigenvector of $K$: $Kv=\lambda v$. Then $\lambda\in \R$ since $K$ is hermitian.
Besides, $K^tJ_{2n}Kv=-Jv$, i.e., $K^t(J_{2n}v)=-1/\lambda J_{2n}v$. By conjugating both sides of this equality we get: $\bar{K}^t(J_{2n}\bar{v})=
-1/\lambda J_{2n}\bar{v}$, i.e., since $K$ is hermitian, $K(J_{2n}\bar{v})=
-1/\lambda J_{2n}\bar{v}$. Therefore if $v$ is an eigenvector of $K$ corresponding to the eigenvalue $\lambda$, then $J_{2n}\bar{v}$ is an eigenvector
of $K$ corresponding to the eigenvalue $-1/\lambda$. Notice that, since $\lambda\in\R$, $\lambda\neq -1/\lambda$.
It follows that one can construct an orthonormal basis $\{v_1, \dots, v_n, J_{2n}\bar{v_1}, \dots, J_{2n}\bar{v_n}\}$ of $\C^{2n}$ consisting of eigenvectors
of $K$. Therefore if we set $W=(v_1\dots v_n -J_{2n}\bar{v}_1\dots -J_{2n}\bar{v}_n)$, $W$ is a unitary matrix such that $K=W\Lambda \overline{W}^t$,
where $\Lambda=diag(\Lambda_1, -\Lambda_1^{-1})$ is a diagonal matrix. One can check that $W$ is symplectic. Then statement $(b)$ follows
using $H=-iK$. 
\hfill$\Box$

\begin{lemma}\label{osp(2,2n)} Consider
 the $3$-bracket  defined by 
(\ref{generalizedC(n)product}) on $M_{1,2n}$ where $H$ is a symplectic matrix. 
\begin{itemize}
\item[(a)] If $H$ is hermitian and $\alpha\in\R^{>0}$ 
(resp.\ $\alpha\in\R^{<0}$) then the 3-algebra $C^3(2n,H;\alpha)$
 is isomorphic to $C^3(2n,H^{2n}_p;1)$
(resp.\ $C^3(2n,H^{2n}_p;-1)$), where $2p$ is the number of positive eigenvalues of $H$.
\item[(b)] If $H$ is anti-hermitian and $\alpha\in i\R^{>0}$ 
(resp.\ $\alpha\in i\R^{<0}$) then 
$C^3(2n,H;\alpha)$ is isomorphic to $C^3(2n,iS^{2n}_n;i)$ 
(resp.\ $C^3(2n,iS^{2n}_n;-i)$).
\end{itemize}
\end{lemma}
{\bf Proof.} 
It
is convenient to identify $M_{1,2n}(\C)$ with the set of matrices of the form
$\left(\begin{array}{cc}
0 & 0\\
\ell & m
\end{array}\right)\in M_{2,2n}$, with $\ell,m\in M_{1,n}$, and
$M_{2n,1}(\C)$ with the set of matrices of the form
$\left(\begin{array}{cc}
n & 0\\
p & 0
\end{array}\right)\in M_{2n,2}$ with $n,p\in M_{n,1}$. Under these 
identifications,
for $Z\in M_{1,2n}(\C)$,
$\psi(Z)=J_{2n}Z^tJ_2^{-1}$.
Let us first assume that $H$ is symplectic and hermitian. Then,
by Lemma \ref{fromFassbenderIkramov} we can write 
$H=yH^{2n}_p\bar{y}^t$ for some matrix $y\in Sp_{2n}(\C)$.
Besides, assume that $\alpha\in\R^{>0}$. Then, for
$h=\left(\begin{array}{cc}
\alpha & 0\\
0 & \alpha^{-1}
\end{array}\right)\in GL_2(\R)$, we can write
$h=x\bar{x}$ for some diagonal matrix $x\in GL_2(\C)$.
Consider the map $\varphi: M_{1,2n}(\C) \rightarrow M_{1,2n}(\C)$ defined by: 
$\varphi(u)=xuy^{-1}$.
We have:
$\varphi([a,b,c]_{H^{2n}_p, 1})=
x[a,b,c]_{H^{2n}_p, 1}y^{-1}=-xaH^{2n}_p\bar{b}^tcy^{-1}+xcH^{2n}_p\bar{b}^tay^{-1}-xc
\psi(a)(\psi(\bar{b}H^{2n}_p))^ty^{-1}$
$=-(xay^{-1})yH^{2n}_p\bar{b}^tx^{-1}(xcy^{-1})+(xcy^{-1})yH^{2n}_p\bar{b}^tx^{-1}(xay^{-1})
-(xcy^{-1})y\psi(a)(\psi(\bar{b}H^{2n}_p))^ty^{-1}=$
$-\varphi(a)y$ $H^{2n}_p\bar{b}^tx^{-1}\varphi(c)+\varphi(c)yH^{2n}_p\bar{b}^tx^{-1}\varphi(a)-
\varphi(c)y\psi(a)(\psi(\bar{b}H^{2n}_p))^t$ $y^{-1}
=-\varphi(a)H\overline{\varphi(b)}^th^{-1}$ $\varphi(c)+
\varphi(c)H\overline{\varphi(b)}^th^{-1}$ $\varphi(a)
-\varphi(c)\psi(\varphi(a))h(\psi(\overline{\varphi(b)}H^t))^t=
[\varphi(a),\varphi(b),$ $\varphi(c)]_{H,\alpha}$\,
since $y\psi(a)(\psi(\overline{b}H^{2n}_p))^ty^{-1}=\psi(\varphi(a))h(\psi(\overline{\varphi(b)}H^t))^t$.
Indeed we have:

\noindent
$\psi(\varphi(a))h(\psi(\overline{\varphi(b)}H^t))^t$ $=
J_{2n}(xay^{-1})^tJ_2^{-1}h(J_{2n}H\overline{\varphi(b)}^tJ_2^{-1})^t$
$=-J_{2n}$ $(y^{-1})^ta^tx^th^{-1}\bar{x}\overline{by^{-1}}H^tJ_{2n}$
$=
-J_{2n}(y^t)^{-1}a^t\bar{b}H^{2n}_py^tJ_{2n}$
$=-yJ_{2n}a^t\bar{b}H^{2n}_pJ_{2n}^ty^{-1}$
$=y\psi(a)(\psi(\bar{b}H^{2n}_p))^ty^{-1}$. Statement $(a)$ with $\alpha\in\R^{>0}$ then follows from Remark \ref{scalars}.
Statement $(a)$ with $\alpha\in\R^{<0}$ can be proven using the same argument with $h=-x\bar{x}$.

One proves statement $(b)$ using the same arguments and the decomposition of an anti-hermitian
symplectic matrix given in Lemma \ref{fromFassbenderIkramov}$(b)$. Indeed, if
$H$ is anti-hermitian, then by Lemma \ref{fromFassbenderIkramov}$(b)$, $H=iyS^{2n}_n\bar{y}^t$ for some symplectic matrix $y$.
If $\alpha\in i\R^{>0}$, i.e., $\alpha=i\beta$ for some $\beta\in\R^{>0}$, we set
$h=\left( 
\begin{array}{cc}
\beta & 0\\
0 & \beta^{-1}
\end{array}\right)$. Then we can write $h=x\bar{x}$ for some diagonal matrix $x\in GL_2(\C)$.
Let $\varphi: M_{1,2n}(\C) \rightarrow M_{1,2n}(\C)$ be defined by $\varphi(u)=xuy^{-1}$. Then,
by arguing as above, one can show that $\varphi([a,b,c]_{iS^{2n}_n,i})=[\varphi(a), \varphi(b), \varphi(c)]_{H, \alpha}$.

A similar argument proves $(b)$ if $\alpha\in i\R^{<0}$.
\hfill$\Box$

\begin{remark}\em
The $N=6$ physical 3-algebra $C^3(2n,H^{2n}_p,1)$ 
can be constructed as explained in Remark \ref{fromalgtophys}, using the following anti-linear involution $C_{n-p}$  of the
algebraic 3-algebra $C^3(2n)$: 
$$C_{n-p}(u)=\bar{u}H^{2n}_{p}, ~u\in M_{1,2n}(\C).$$
Therefore, we have $C^3(2n,H^{2n}_p,1)=C^3(2n)_{ph, C_{n-p}}$. 
\end{remark}


%

\begin{example}\label{K(m,2)}\em
%
Consider the generalised Poisson algebra $P(m,0)$  in the (even) 
indeterminates $p_1, \dots,$ $p_k, q_1, \dots, q_k$
(resp.\ $p_1, \dots,p_k, q_1, \dots, q_k, t$) if $m=2k$ (resp.\ $m=2k+1$)
endowed with the bracket: 
\begin{equation} 
\{f,g\}=(2-E)(f)\frac{\partial g}{\partial t}-\frac{\partial f}{\partial t}(2-E)(g)+\sum_{i=1}^k(\frac{\partial f}{\partial p_i}
\frac{\partial g}{\partial q_i}-
\frac{\partial f}{\partial q_i}\frac{\partial g}{\partial p_i}),
\label{Poisson}
\end{equation}
where $E=\sum_{i=1}^k(p_i\frac{\partial}{\partial p_i}+q_i\frac{\partial}
{\partial q_i})$
(the first two terms in (\ref{Poisson}) vanish if $m$ is even) and the 
derivation $D=2\frac{\partial}{\partial t}$ (which is 0 if $m$ is even).
Let
$\sigma_{\varphi}: f(p_i,q_i)\mapsto -{f}(\varphi (p_i),\varphi(q_i))$ 
(resp. $\sigma_{\varphi}: f(t,p_i,q_i)\mapsto -{f}(\varphi(t),\varphi (p_i),\varphi(q_i))$, 
where $\varphi$ is an involutive linear change of variables 
(i.e. $\varphi^2=1$),
multiplying by $-1$ the 1-form  
$\sum_i (p_i dq_i - q_i dp_i)$ if $m$ is even 
(resp. $dt + \sum_i (p_i dq_i - q_i dp_i)$ if $m$ is odd).
Then the $3$-bracket 
\begin{equation}
[f,g,h]=\{f, \sigma_{\varphi}(g)\}h+\{f,h\}\sigma_{\varphi}(g)+f\{\sigma_{\varphi}(g),h\}+
D(f)\sigma_{\varphi}(g)h-f\sigma_{\varphi}(g)D(h),
\label{eq:1.1}
\end{equation}
defines on $P(m,0)$ an algebraic $N=6$ $3$-algebra structure that we denote by $P^3(m,\varphi)$ \cite[Example 1.6]{CantaK6}.

If in (\ref{eq:1.1}) we replace $\sigma_{\varphi}$ with the map $\bar{\sigma}_{\varphi}$
 defined by:
$\bar{\sigma}_{\varphi}: f(p_i,q_i)\mapsto -\bar{f}(\varphi (p_i),\varphi(q_i))$ 
(resp. $\bar{\sigma}_{\varphi}: f(t,p_i,q_i)\mapsto -\bar{f}(\varphi(t),\varphi (p_i),\varphi(q_i))$, 
we get a physical $N=6$ $3$-algebra that we denote by $P^3(m,\varphi;\bar{\,})_{+}$.
Similarly, we shall denote by $P^3(m,\varphi;\bar{\,})_{-}$ the physical $N=6$ $3$-algebra
with bracket (\ref{eq:1.1}) where $\sigma_{\varphi}$ is replaced with the map $-\bar{\sigma}_{\varphi}$.

Note that $\overline{\sigma_{\varphi}(g)}=\sigma_{\varphi}(\bar{g})$ if and only if $\varphi$ is a real change of variables, 
hence in this case the complex conjugation $C_0: P^3(m,\varphi)\rightarrow P^3(m,\varphi)$, $C(f)=\bar{f}$, is an anti-linear involution
of the algebraic $N=6$ 3-algebra $P^3(m,\varphi)$ and, by Remark \ref{fromalgtophys},  $P^3(m,\varphi;\bar{\,})_{\pm}=P^3(m,\varphi)_{ph,\pm C_0}$. 
\end{example}

%

\begin{example}\label{S(1,2)new}\em
Let $A=\C[[x]]^{\langle 1\rangle}\oplus\C[[x]]^{\langle 2\rangle}$ be the 
direct sum of two copies of the algebra $\C[[x]]$; for $f\in\C[[x]]$, 
denote by $f^{\langle i\rangle}$ the
corresponding element in $\C[[x]]^{\langle i\rangle}$. Let $D=d/dx$, and let 
$a=(a_{ij})\in M_{2,2}(\C)$.

The following $3$-bracket defines on $A$ an algebraic $N=6$ 3-algebra structure for every
$a\in SL_2(\C)$ such that either $a^2=1$ (and in this case $\varphi=1$) or $a^2=-1$ (and in this case $\varphi=-1$):
for $i,j=1$ or $2$,
$$[f^{\langle i\rangle},g^{\langle i\rangle},h^{\langle i\rangle}]=
(-1)^ia_{ij}((fD(h)-D(f)h)g(\varphi(x)))^{\langle i\rangle}
~~{\mbox{for}}~~j\neq i;$$
$$[f^{\langle i\rangle},g^{\langle j\rangle},h^{\langle i\rangle}]=
(-1)^ia_{jj}((fD(h)-D(f)h)g(\varphi(x)))^{\langle i\rangle}
~~{\mbox{for}}~~j\neq i;$$
$$[f^{\langle 1\rangle},g^{\langle j\rangle},h^{\langle 2\rangle}]=
a_{j1}((fD(g(\varphi(x)))-D(f)g(\varphi(x)))h)^{\langle 1\rangle})+
a_{j2}(f(hD(g(\varphi(x)))-D(h)g(\varphi(x))))^{\langle 2\rangle}.$$
This bracket is extended to $A$ by skew-symmetry in the first and third entries.
We denote this 3-algebra by $SW^3(a)$ \cite[Example 1.7]{CantaK6}. 

Similarly, for every $a\in SL_2(\C)$ such that either $\bar{a}a=1$ (and in this case $\lambda=\pm 1$) or $\bar{a}a=-1$ (and in this case $\lambda=\pm i$): we can define a physical $N=6$ 3-algebra structure by setting
$$[f^{\langle i\rangle},g^{\langle i\rangle},h^{\langle i\rangle}]=\lambda(\exp(-t))
(-1)^i\overline{a_{ij}}((fD(h)-D(f)h)\overline{g(\varphi_t(x))})^{\langle i\rangle}
~~{\mbox{for}}~~j\neq i;$$
$$[f^{\langle i\rangle},g^{\langle j\rangle},h^{\langle i\rangle}]=\lambda(\exp(-t))
(-1)^i\overline{a_{jj}}((fD(h)-D(f)h)\overline{(g(\varphi_t(x))})^{\langle i\rangle}
~~{\mbox{for}}~~j\neq i;$$
$$[f^{\langle 1\rangle},g^{\langle j\rangle},h^{\langle 2\rangle}]=\lambda(\exp(-t))
(\overline{a_{j1}}((f\overline{D(g(\varphi_t(x)))}-D(f)\overline{g(\varphi_t(x))})h)^{\langle 1\rangle})$$
$$\,\,\,\,\,\,\,\,\,\,\,+\overline{a_{j2}}(f(h\overline{D(g(\varphi_t(x)))}-D(h)\overline{g(\varphi_t(x))}))^{\langle 2\rangle}),$$
where $t\in i\R$ and  $\varphi_t(x)=\exp(2t)x$,
and extending it to $A$ by skew-symmetry in the first and third entries.
We shall denote this physical $N=6$ 3-algebra by $SW^3(a;\varphi_t)_{\pm}$, depending on $\lambda$ being equal to $\pm 1$ (resp.\ $\pm i$)
if $\bar{a}a=1$ (resp.\ $\bar{a}a=-1$). 

Note that if $C_0$ denotes the complex conjugation of $\C[[x]]$, i.e.,
for $f\in\C[[x]]$, $f=\sum_{i\in\Z}\alpha_ix^i$,
$C_0(f)= \sum_{i\in\Z}\overline{\alpha}_ix^i$, then $C_0$ is an anti-linear involution of
$SW^3(a)$ if and only if $a\in SL_2(\R)$. Hence, by Remark \ref{fromalgtophys}, if $a\in SL_2(\R)$ and $t=k\pi i$, $k\in\Z$, then $SW^3(a;\varphi_t)_{\pm}=SW^3(a)_{ph,\pm C_0}$.

\end{example}

\begin{example}\label{SKO(2,3;1)new}\em
Let $A=\C[[x_1,x_2]]$ and $D_i= \frac{\partial}{\partial x_i}$ for $i=1,2$.
Consider the following $3$-bracket:
\begin{equation}
[f,g,h]=\det
\left(
\begin{array}{ccc}
f & \varphi(g) & h\\
D_1(f) & D_1(\varphi(g)) & D_1(h)\\
D_2(f) & D_2(\varphi(g)) & D_2(h)
\end{array}\right),
\label{scircsign}
\end{equation}
where 
$\varphi$ is an  automorphism of the associative algebra $A$.
If $\varphi$ is a linear change of variables with determinant equal to 1, and 
$\varphi^2=1$, then $A$ with $3$-bracket (\ref{scircsign}) is an algebraic
$N=6$ 
$3$-algebra, which we denote by $W^3(\varphi)$ \cite[Example 1.8]{CantaK6}.  
Likewise, consider the following $3$-bracket on $A$:
\begin{equation}
[f,g,h]=\det
\left(
\begin{array}{ccc}
f & \overline{\varphi(g)} & h\\
D_1(f) & D_1(\overline{\varphi(g)}) & D_1(h)\\
D_2(f) & D_2(\overline{\varphi(g)}) & D_2(h)
\end{array}\right)
\label{scircsignph}
\end{equation}
where $\varphi$ is an automorphism of the algebra $A$.
If $\varphi$ is a linear change of variables 
such that
$\bar{\varphi}\varphi=1$ 
then $A$ with the bracket (\ref{scircsignph}) is a physical  $N=6$ 
$3$-algebra that we shall denote by $W^3(\varphi, \bar{\,})_{+}$. Likewise we shall denote by $W^3(\varphi, \bar{\,})_{-}$
the $3$-algebra $A$ with the negative of the 3-bracket (\ref{scircsignph}).

Note that $\overline{\varphi(g)}=\varphi({\bar{g}})$ if and only if $\varphi$ is a real change of variables, hence in this case the complex
conjugation of $A$ is an anti-linear involution of the algebraic $N=6$ 3-algebra $W^3(\varphi)$, and, by Remark \ref{fromalgtophys},
$W^3(\varphi, \bar{\,})_{\pm}=W^3(\varphi)_{ph, \pm C_0}$.
\end{example}

\begin{example}\label{SKO(2,3;beta)new}\em
Let $A=\C[[x_1,x_2]]$ and $D_i= \frac{\partial}{\partial x_i}$ for $i=1,2$.
Consider the following $3$-bracket:
\begin{equation}
[f,g,h]=\det
\left(
\begin{array}{ccc}
(2-E)(f) & (2\beta-E)(\overline{\varphi(g)}) & (2-E)(h)\\
D_1(f) & D_1(\overline{\varphi(g)}) & D_1(h)\\
D_2(f) & D_2(\overline{\varphi(g)}) & D_2(h)
\end{array}\right)
\label{betanambu}
\end{equation}
where $\beta\in\C$, $E=\sum_{i=1}^2x_i\frac{\partial}{\partial x_i}$ and $\varphi$ is a linear change of indeterminates.
If $\varphi\in GL_2$ is such that $\bar{\varphi}\varphi=I_2$, $\beta\notin\R$ and $|\beta|=1$, then $A$ with $3$-bracket (\ref{betanambu}) is a physical
$N=6$ 
$3$-algebra that we shall denote by $W^3_{\beta}(\varphi)_+$.  Likewise, we shall denote by $W^3_{\beta}(\varphi)_-$ the $3$-algebra $A$ with the negative
of the 3-bracket (\ref{betanambu}). 
\end{example}

\begin{example}\label{SHO(3,3)new}\em
Let $A=\C[[x_1,x_2,x_3]]$ and $D_i=\frac{\partial}{\partial x_i}$ for 
$i=1,2,3$.  Consider the following $3$-bracket on $A$:
\begin{equation}
[f,g,h]=\det
\left(
\begin{array}{ccc}
D_1(f) & D_1(\varphi(g)) & D_1(h)\\
D_2(f) & D_2(\varphi(g)) & D_2(h)\\
D_3(f) & D_3(\varphi(g)) & D_3(h)
\end{array}\right)
\label{scircfi}
\end{equation}
where $\varphi$ is an automorphism of the algebra $A$.
If $\varphi$ is a linear change of variables with determinant equal to 1, and 
$\varphi^2=1$, then $A$ with the bracket (\ref{scircfi}) is an algebraic  $N=6$ 
$3$-algebra and $\C 1$ is an ideal of this 3-algebra \cite[Example 1.9]{CantaK6}. We denote by $S^3(\varphi)$ the quotient
3-algebra
$A/\C 1$.  

Likewise, consider the following $3$-bracket on $A$:
\begin{equation}
[f,g,h]=\alpha\det
\left(
\begin{array}{ccc}
D_1(f) & D_1(\overline{\varphi(g)}) & D_1(h)\\
D_2(f) & D_2(\overline{\varphi(g)}) & D_2(h)\\
D_3(f) & D_3(\overline{\varphi(g)}) & D_3(h)
\end{array}\right)
\label{scircfiph}
\end{equation}
where $\alpha\in\C$ and $\varphi$ is an automorphism of the associative algebra $A$. 
If $\varphi$ is a linear change of variables with determinant equal to $1$ (resp.\ $-1$),
$\alpha=\pm 1$ (resp.\ $\alpha=\pm i$) and 
$\bar{\varphi}\varphi=1$, then $A$ with the bracket (\ref{scircfiph}) is a physical  $N=6$ 
$3$-algebra and $\C 1$ is an ideal of this 3-algebra. We shall denote by $S^3(\varphi, \bar{\,})_{\pm}$ the quotient
3-algebra $A/\C 1$.  

Note that $\overline{\varphi(g)}=\varphi({\bar{g}})$ if and only if $\varphi$ is a real change of variables, hence in this case the complex
conjugation of $A$ is an anti-linear involution of the algebraic $N=6$ 3-algebra $S^3(\varphi)$, and, by Remark \ref{fromalgtophys},
$S^3(\varphi, \bar{\,})_{\pm}=S^3(\varphi)_{ph, \pm C_0}$.
\end{example}

The main result of the paper is the following theorem.
\begin{theorem}\label{main}
The following is a complete list, up to isomorphisms, of simple linearly compact physical $N=6$
$3$-algebras:
\begin{itemize}
\item[-] finite-dimensional:
\begin{itemize}
\item[(i)] $A^3(m,n;t)_{ph,C_{p,q}}$ ($0\leq p\leq m$, $0\leq q\leq n$); 
\item[(ii)] $A^3(n)_{\pm}$;
\item[(iii)] $C^3(2n)_{ph, \pm C_p}$ ($1\leq p\leq n$); 
\item[(iv)] $C^3(2n, iS^{2n}_n, \pm i)$.
\end{itemize}
\item[-] infinite-dimensional:
\begin{itemize}
\item[(a)] $P^3(m, \varphi;\bar{\,})_{\pm}$ ($m \geq 1$);
\item[(b)] $SW^3(a, \varphi_t)_{\pm}$;
\item[(c)] $W^3(\varphi, \bar{\,})_{\pm}$; 
\item[(d)] $S^3(\varphi, \bar{\,})_{\pm}$;
\item[(e)] $W^3_{\beta} (\varphi)_\pm$.
\end{itemize}
\end{itemize}

Among these a complete list of physical $N=6$ 3-algebras which are trivially related to algebraic $N=6$ 3-algebras over $\C$, is the following:
\begin{itemize}\item[-] finite-dimensional:
\begin{itemize}
\item[(i)]
$A^3(m,n;t)_{ph,C_{p,q}}$ ($0\leq p\leq m$, $0\leq q\leq n$);
\item[(ii)] $C^3(2n)_{ph, \pm C_p}$ ($1\leq p\leq n$).
\end{itemize}
\item[-] infinite-dimensional:
\begin{itemize}
\item[($a^\prime$)] 
$P^3(m, \varphi;\bar{\,})_{\pm}=P^3(m,\varphi)_{ph,\pm C_0}$ ($m \geq 1$) where $\varphi$ is a real change of variables;
\item[($b^\prime$)] $SW^3(a, \varphi_t)_{\pm}=SW^3(a)_{ph,\pm C_0}$ where $a\in SL_2(\R)$ and $t=k\pi i$, $k\in\Z$;
\item[($c^\prime$)] $W^3(\varphi, \bar{\,})_{\pm}=W^3(\varphi)_{ph, \pm C_0}$ where $\varphi$ is a real change of variables
\item[($d^\prime$)] $S^3(\varphi, \bar{\,})_{\pm}=S^3(\varphi)_{ph, \pm C_0}$ where $\varphi$ is a real change of variables.
\end{itemize}
\end{itemize}
\end{theorem}
{\bf Proof.}
Theorem \ref{th:2.3} from Section 2 reduces the classification in question to that of 
the pairs $(L,\sigma)$, where
$L=L_{-1}\oplus L_0 \oplus L_1$
is a simple linearly compact Lie superalgebra with a consistent $\ZZ$-grading satisfying properties $(i)$ and $(ii)$
of Theorem \ref{th:2.3},
and $\sigma$ is an anti-linear graded conjugation of $L$. A complete list of possible such
Lie superalgebras
$L=L_{-1}\oplus L_0 \oplus L_1$ is given by Proposition \ref{fdshortconsistent} from Section \ref{section3}.
Finally, a complete list of anti-linear graded conjugations of these $L$ is given, in the finite-dimensional case by
Proposition \ref{gradedanticonjforfds} from Section \ref{section3}, and in the infinite-dimensional case by Proposition
\ref{listofgradedconj} from Section \ref{section3}.

By Theorem \ref{th:2.3}$(b)$, the physical $N=6$ $3$-algebra is identified with $\Pi L_{-1}$,
on which the $3$-bracket is given by the formula 
$[a,b,c]=[[a,\sigma(b)],c]$. This formula, applied to the $\ZZ$-graded finite-dimensional
Lie superalgebras $L$ with graded conjugations, described by Proposition \ref{gradedanticonjforfds}
$(a), (b)$ and $(c)$ produces the $3$-algebras
$A^3(m,n;t)_{ph, C_{p,q}}$, $A^3(n)_{\pm}$ and $C^3(2n)_{ph, \pm C_p}$, $C^3(2n,iS^{2n}_n; \pm i)$.
The same formula, applied to the $\ZZ$-graded infinite-dimensional
Lie superalgebras $L$ with graded conjugations, described by Proposition \ref{listofgradedconj}
$(a), (b), (c), (d), (e)$ and $(f)$ produces the $3$-algebras
$P^3(2k, \varphi;\bar{\,})_{\pm}$ ($m \geq 1$), $P^3(2k+1, \varphi;\bar{\,})_{\pm}$,
$SW^3(a, \varphi_t)_{\pm}$,
$W^3(\varphi, \bar{\,})_{\pm}$, $S^3(\varphi, \bar{\,})_{\pm}$, $W^3_{\beta}(\varphi)_\pm$.

 The fact that all of them are indeed phyisical $N=6$ $3$-algebras follows
automatically from Theorem \ref{th:2.3}$(b)$. 
\hfill$\Box$

\begin{remark}\em As we already noticed at the end of Example \ref{sl(2,2)}, we have $A^3(m,n;t)_{ph,C_{0,0}}=A^3(m,n;t)_{ph,C_{n,m}}$
and $A^3(m,n;t)_{ph,C_{n,0}}=A^3(m,n;t)_{ph,C_{0,m}}$.
\end{remark}


\section{Palmkvist's construction}
This section is based on the ideas of \cite{P}.
\begin{definition}
An element $b$ of an $N=6$ 3-algebra $\g$ is called central if $[a, b, c]=0$ for all $a$ and $c$ in $\g$; the
subspace consisting of central elements is called the center of $\g$.
\end{definition}
\begin{remark}\em
It is immediate to see from the
axioms that the center is an ideal of $\g$. In particular the center of a simple $N=6$ 3-algebra is zero
(by definition the one-dimensional 3-algebra with zero 3-bracket is not simple). 
\end{remark}
\begin{definition}
Let $\g=\oplus_{j\in\Z} \g_j$ be a 
Lie superalgebra over $\C$ with a consistent $\Z$-grading. A graded conjugation (resp.\ anti-linear graded conjugation) 
of $\g$ is a linear (resp.\ anti-linear) Lie
superalgebra automorphism
$\varphi: \g \rightarrow \g$ such that
\begin{enumerate}
\item $\varphi(\g_j)=\g_{-j}$
\item $\varphi^2(x)=(-1)^kx ~{\mbox{for}}~ x\in\g_k$.
\end{enumerate}
\end{definition}
\begin{theorem}\label{tel}
Let $\g=\g_{-1}\oplus\g_0\oplus\g_1$ be a consistently $\Z$-graded 
Lie superalgebra with an anti-linear graded conjugation  $\varphi$. 
Then the $3$-bracket
\begin{equation}
[u,v,w]:=[[u,\varphi(v)],w]
\label{palm}
\end{equation}
defines on $\Pi\g_{-1}$ a physical $N=6$ 3-algebra structure (here and further $\Pi$ stands for
the parity reversal).
\end{theorem}
{\bf Proof.} Bracket (\ref{palm}) is obviously linear in the first and third argument and anti-linear in the second one.
Since the grading of $\g$ is consistent, $\g_{-1}$ and $\g_1$ are completely odd and $\g_0$
is even. For $u,v,w\in \g_{-1}$ we thus have:
$[u,v,w]:=[[u,\varphi(v)],w]=[u,[\varphi(v),w]]=-[[\varphi(v),w],u]=-[[w,\varphi(v)],u]=-[w,v,u]$,
which proves anti-commutativity.

Besides, for $u,v,x,y,z\in\g_{-1}$ we have:
$[u,v,[x,y,z]]-[x,y,[u,v,z]]=[[u,\varphi(v)],[[x,\varphi(y)],z]]-[[x,\varphi(y)],[[u,\varphi(v)],z]]=$
$[[u,\varphi(v)],[x,\varphi(y)]],z]+[[x,\varphi(y)],[[u,\varphi(v)],z]]-[[x,\varphi(y)],[[u,\varphi(v)],z]]=$
$[[u,\varphi(v)],[x,\varphi(y)]],z]=[[[u,\varphi(v)],x],\varphi(y)],z]+
[[[x,[[u,\varphi(v)],\varphi(y)]],z]=[[[u,\varphi(v)],x],\varphi(y)],z]-
[[[x,\varphi([[\varphi(u),v],y])],z]=$
$[[[u,v,x],y,z]-[x,[v,u,y],z]$
\hfill$\Box$

\bigskip

We shall now associate to a physical $N=6$ 3-algebra
$T$ with zero center a $\Z$-graded Lie superalgebra 
$Lie T=Lie_{-1} T\oplus Lie_0 T\oplus Lie_1 T$, as follows. For $x,y\in T$, 
denote by $L_{x,y}$ the endomorphism
of $T$ defined by $L_{x,y}(z)=[x,y,z]$. 
Besides, for $x\in T$, denote by $\varphi_x$ the map in 
$Hom(\Pi T\otimes \Pi T, \Pi T)$
defined by $\varphi_x(y,z)=-[y,x,z]$. 

We let $Lie_{-1} T=\Pi T$, $Lie_0 T=\langle L_{x,y} ~|~ x,y\in T\rangle$, 
$Lie_1 T=\langle \varphi_x ~|~ x\in T\rangle$,
and let $Lie T=Lie_{-1}T \oplus Lie_0 T \oplus Lie_1T$.  
Define the map $\sigma: Lie T \longrightarrow Lie T$ by setting
$(x,y,z\in T)$:
$$z\mapsto -\varphi_z,
~~~\varphi_z\mapsto z, \,\,\,L_{x,y}\mapsto -L_{y,x},$$
and extending it on $Lie T$ by anti-linearity.
The map $\sigma$ is well defined since the center of $T$ is zero. Indeed, 
by definition,
 $\varphi_z(x,y)=-[x,z,y]$, hence $\varphi_z=0$ only for $z=0$.
Besides, $L_{a,b}=0$ for some $a,b\in T$ implies, by the fundamental identity, 
that $[b,a,y]$ is a central element for any $y$, hence
$L_{b,a}=0$. Besides, the following relations hold:
 $\varphi_{\alpha z}=\overline{\alpha}\varphi_z$ and $\alpha L_{x,y}=L_{\alpha x,y}=L_{x,\overline{\alpha}y}$, which are
 consistent with the anti-linearity of $\sigma$.
 
\begin{theorem}
\label{th:2.3}
(a) $Lie T$ is a $\Z$-graded Lie superalgebra with a short consistent grading,
satisfying the following two properties:

(i) any non-zero $\ZZ$-graded ideal of $Lie T$ has a non-zero 
intersection with both $Lie_{-1}T$ and $Lie_1 T$;

(ii) $[Lie_{-1} T, Lie_1 T]=Lie_0 T$.

(b) $\sigma$ is an anti-linear graded conjugation of the $\ZZ$-graded Lie superalgebra
$Lie T$ and the 3-bracket on $T$ is recovered from the bracket on $Lie T$
by the formula $[x,y,z]=[[x,\sigma(y)],z]$.

(c) The correspondence $T\longrightarrow (Lie T, \sigma)$ 
is bijective between the isomorphism classes of physical
$N=6$ $3$-algebras with zero center and the isomorphism classes of the pairs $(Lie T,\sigma)$, where  
$Lie T$ is a $\Z$-graded Lie superalgebra with a short consistent grading,
satisfying properties (i) and (ii), and $\sigma$ is an anti-linear graded conjugation
of $Lie T$.

(d) A $3$-algebra $T$ is simple (resp. finite-dimensional or linearly compact)
if and only if $Lie T$ is.  
\end{theorem}
{\bf Proof.} For $x,y,z\in T$, we have:
$[\varphi_x, z]=-L_{z,x}$, and
\begin{equation}
[L_{x,y},L_{x',y'}]=L_{[x,y,x'],y'}-L_{x',[y,x,y']}.
\label{eq:1.2}
\end{equation}
Note that
$[L_{x,y},\varphi_z]=-\varphi_{[y,x,z]}$.
Finally, $[[\varphi_x, \varphi_y], z]=[\varphi_x, [\varphi_y,z]]+[\varphi_y, [\varphi_x,z]]=
-[\varphi_x, L_{z,y}]-[\varphi_y, L_{z,x}]=[L_{z,y}, \varphi_x]+[L_{z,x}, \varphi_y]=
-\varphi_{[y,z,x]}-\varphi_{[x,z,y]}=0$. Hence $[\varphi_x, \varphi_y]=0$.

It follows that $Lie T$ is indeed a $\ZZ$-graded Lie superalgebra,
satisfying $(ii)$ and such that any non-zero ideal has a non-zero
intersection with $Lie_{-1}T$. It is straightforward to check $(b)$,
hence any non-zero ideal of $Lie T$ has a non-zero intersection
with $Lie_1T$, which completes the proof of $(a)$.

(c) follows from Theorem \ref{tel} where the inverse of the correspondence $T\longrightarrow (Lie T, \sigma)$
is given. Finally, since the simplicity of $T$, by definition,
means that all operators $L_{x,y}$ have no common non-trivial invariant 
subspace in $T$, it follows that $T$ is a simple $3$-algebra if and only if
$Lie_0 T$ acts irreducibly on $Lie_{-1} T$. Hence, by the properties $(i)$ and 
$(ii)$ of $Lie T$, $T$ is simple if and only if $Lie T$ is simple. The rest of 
(d) is clear as well. 
\hfill$\Box$


\section{Classification of anti-linear graded conjugations}\label{section3}
In this section we shall classify all anti-linear graded conjugations $\sigma$ of all 
$\Z$-graded simple linearly compact Lie superalgebras $\g$ with a short 
consistent grading
$\g=\g_{-1}\oplus\g_0\oplus\g_1$. 

In \cite[\S 3]{CantaK6} we classified all short gradings of all simple linearly
compact Lie superalgebras and deduced the following lists of
simple finite-dimensional $\Z$-graded Lie superalgebras with a short consistent grading
 $\g=\g_{-1}\oplus \g_0\oplus \g_1$
such that $\g_{-1}$ and $\g_1$ have the same dimension
 \cite[Remark 3.2]{CantaK6} and 
of simple infinite-dimensional $\Z$-graded Lie superalgebras with a short consistent grading
 $\g=\g_{-1}\oplus \g_0\oplus \g_1$
such that $\g_{-1}$ and $\g_1$ have the same growth and size
\cite[Remark 3.4]{CantaK6}:
\begin{proposition}\label{fdshortconsistent} \begin{itemize}
\item[(a)] A complete list of simple finite-dimensional $\Z$-graded Lie superalgebras with a short consistent grading
 $\g=\g_{-1}\oplus \g_0\oplus \g_1$
such that $\g_{-1}$ and $\g_1$ have the same dimension, is, up
to isomorphism, as follows:
\begin{itemize}
\item[-] $psl(m,n)$ with $m,n\geq 1$, $m+n\geq 2$,
with the grading $f(\epsilon_1)=\dots=f(\epsilon_m)=1$,
$f(\delta_1)=\dots=f(\delta_n)=0$;
\item[-] $osp(2,2n)$, $n\geq 1$, with the grading $f(\delta_i)=0$ for all $i$,
$f(\epsilon_1)=1$.
\end{itemize}
\item[(b)] 
A complete list of simple infinite-dimensional $\Z$-graded Lie superalgebras with a short consistent grading
 $\g=\g_{-1}\oplus \g_0\oplus \g_1$
such that $\g_{-1}$ and $\g_1$ have the same growth and size, is, up
to isomorphism, as follows:
\begin{itemize}
\item[-] $S(1,2)$ with the grading of type $(0|1,1)$;
\item[-] $H(2k,2)$ with the grading of type $(0,\dots,0|1,-1)$;
\item[-] $K(2k+1,2)$ with the grading of type $(0,\dots,0|1,-1)$;
\item[-] $SHO(3,3)$ with the grading of type $(0,0,0|1,1,1)$;
\item[-] $SKO(2,3;\beta)$ with the grading of type $(0,0|1,1,1)$.
\end{itemize}
\end{itemize}
\end{proposition}
In Examples \ref{gcforsl(n,n)} and \ref{gradedconjforfds} the $\Z$-grading of $\g$ is given by
Proposition \ref{fdshortconsistent}(a). 
\begin{example}\em\label{gcforsl(n,n)}
Let $\g=psl(n,n)$. Then the maps $\tau_{\pm}\left(
\begin{array}{cc}
a & b\\
c & d
\end{array}\right)=\left(
\begin{array}{cc}
\bar{d} & \pm i\bar{c}\\
\mp i\bar{b} & \bar{a}
\end{array}\right)$
are anti-linear graded conjugations of $\g$.
\end{example}

\begin{remark}\label{Ccomposed}\em
Let $\sigma$ be a graded conjugation of $\g$ and let $C$ be an anti-linear involution
of $\g=\g_{-1}\oplus\g_0\oplus\g_1$ preserving the grading and such that
$\sigma\circ C=C\circ\sigma$. Then $\tilde{\sigma}:=C\circ\sigma$ is an anti-linear
graded conjugation of $\g$.
\end{remark}

%
%
\begin{example}\label{gradedconjforfds}\em
The following examples of anti-linear graded conjugations can be deduced from \cite[Propositions 4.3, 4.8]{CantaK6}, using
Remark \ref{Ccomposed}:
\begin{itemize}
\item[(a)] $\g=psl(m,n)$: $\tilde{\sigma}_1\left(
\begin{array}{cc}
a & b\\
c & d
\end{array}\right)=\left(
\begin{array}{cc}
-\bar{a}^t & \bar{c}^t\\
-\bar{b}^t & -\bar{d}^t
\end{array}\right)$.
\item[(b)] $\g=osp(2,2n)$: $\tilde{\sigma}_1$.
\end{itemize}
\end{example}

\begin{remark}\label{howtogetall}\em
Let $\varphi_1$ and $\varphi_2$ be two anti-linear graded conjugations of a consistently $\Z$-graded Lie superalgebra
$\g=\g_{-1}\oplus\g_0\oplus\g_1$. Let $f=\varphi_2\varphi_1^{-1}$. Then $f$ is a grading preserving linear automorphism
of $\g$ such that 
\begin{equation}
(f\circ\varphi_1)^2(x_j)=(-1)^j x_j ~~{\mbox{for}}~x_j\in\g_j.
\label{automorphismproperty}
\end{equation} It follows that if $\varphi_1$ is an anti-linear
graded conjugation of $\g$, then any other such conjugation is of the form $f\circ \varphi_1$ for some grading preserving
linear automorphism $f$ of $\g$ satisfying (\ref{automorphismproperty}).  
\end{remark}

\begin{definition}\label{equivalent}
Two (linear or anti-linear) graded conjugations $\sigma_1$ and $\sigma_2$ of a Lie superalgebra $\g=\g_{-1}\oplus\g_0\oplus \g_1$ are
called equivalent if $\sigma_2=\varphi \sigma_1 \varphi^{-1}$, where $\varphi$ is a  grading
preserving automorphism of $\g$.
\end{definition}

\begin{proposition}\label{gradedanticonjforfds}
The following is a complete list, up to equivalence, of anti-linear graded conjugations of all simple finite-dimensional Lie superalgebras with the $\Z$-gradings from Proposition \ref{fdshortconsistent}(a): 
\begin{itemize}
\item[(a)] $\g=psl(m,n)$: $Ad\,diag(S^m_p,S^n_q)\circ\tilde{\sigma}_1$
for $p=0,\dots, m$ and $q=0,\dots, n$, where $S^m_p$ and $S^n_q$ are defined in Example \ref{sl(2,2)}.
\item[(b)] $\g=psl(n,n)$: $\sigma=\tau_{\pm}$;
\item[(c)] $\g=osp(2,2n)$: $Ad\,diag(\pm A,H)\circ\tilde{\sigma}_1$ where either $A=I_2$ and
$H=H^{2n}_p$ for some $p=0,\dots,n$, or $A=diag(i,-i)$ and $H=iS^{2n}_n$.
\end{itemize}\end{proposition}
{\bf Proof.}
The Lie superalgebra $sl(m,n)$ has a short consistent grading such that 
$\g_0=\g_{\0}$ consists of
matrices of the form 
$\left(
\begin{array}{cc}
\alpha & 0\\
0 & \delta
\end{array}\right)$, where  
$tr \alpha=tr\delta$, $\g_{-1}$ is the set of matrices of the form
$\left(
\begin{array}{cc}
0 & 0\\
\gamma & 0
\end{array}\right)$, and
$\g_{1}$ is the set of matrices of the form
$\left(
\begin{array}{cc}
0 & \beta\\
0 & 0
\end{array}\right)$.

For $m\neq n$ every automorphism of $\g=sl(m,n)$ is either of the form $Ad~diag(A, B)$ for some matrices $A\in GL_m(\C)$, $B\in GL_n(\C)$,
or of the form $Ad~diag(A, B)\circ \sigma_1$ \cite{S}, where, for
$\left(
\begin{array}{cc}
a & b\\
c & d
\end{array}\right)\in sl(m,n)$,
$\sigma_1\left(
\begin{array}{cc}
a & b\\
c & d
\end{array}\right)=\left(
\begin{array}{cc}
-{a}^t & {c}^t\\
-{b}^t & -{d}^t
\end{array}\right).$
Note that
$Ad~diag(A,B)\left(
\begin{array}{cc}
a & b\\
c & d
\end{array}\right)=\left(
\begin{array}{cc}
AaA^{-1} & AbB^{-1}\\
BcA^{-1} & BdB^{-1}
\end{array}\right)$, hence
$Ad~diag(A,B)\circ\sigma_1\left(
\begin{array}{cc}
a & b\\
c & d
\end{array}\right)=\left(
\begin{array}{cc}
-A{a}^tA^{-1} & A{c}^tB^{-1}\\
-B{b}^tA^{-1} & -B{d}^tB^{-1}
\end{array}\right).$
Therefore only automorphisms of the form $Ad~diag(A,B)$ preserve the short consistent grading of $\g$.
By Remark \ref{howtogetall} and Example \ref{gradedconjforfds}, every anti-linear graded conjugation of
$\g$ is of the form $\varphi=Ad~diag(A,B)\circ\tilde{\sigma}_1$. Condition (\ref{automorphismproperty})
then implies $A=\lambda \bar{A}^t$ and $B=\lambda\bar{B}^t$, for some $\lambda\in\C$ such that
$|\lambda|=1$. Indeed we have:
$\varphi\left(
\begin{array}{cc}
a & b\\
c & d
\end{array}\right)=\left(
\begin{array}{cc}
-A\bar{a}^tA^{-1} & A\bar{c}^tB^{-1}\\
-B\bar{b}^tA^{-1} & -B\bar{d}^tB^{-1}
\end{array}\right),$ hence
$\varphi^2\left(
\begin{array}{cc}
a & 0\\
0 & d
\end{array}\right)=\left(
\begin{array}{cc}
A(\bar{A}^{-1})^ta(\bar{A})^tA^{-1} & 0\\
0 & B(\bar{B}^{-1})^td(\bar{B})^tB^{-1}
\end{array}\right)$
from which we get $A=\lambda\bar{A}^t$ and $B=\mu \bar{B}^t$ for some $\lambda, \mu\in\C$. 
Besides,
$\varphi^2\left(
\begin{array}{cc}
0 & b\\
0 & 0
\end{array}\right)=\left(
\begin{array}{cc}
0 & -A(\bar{A}^{-1})^tb\bar{B}^tB^{-1}\\
0 & 0
\end{array}\right),$
hence $\lambda=\sigma$. Finally, taking the conjugate transpose of both sides of the equality $A=\lambda\bar{A}^t$, we get
$\bar{A}^t=\bar{\lambda}A=\bar{\lambda}\lambda\bar{A}^t$, i.e., $|\lambda|=1$. Statement $(a)$ then follows from
Lemma \ref{sl(m,n)}.

If $m=n$, in addition to the automorphisms described above, $psl(n,n)$ has automorphisms of the form $Ad~diag(A, B)\circ \Pi$,
$Ad~diag(A, B)\circ \Pi\circ\sigma_1$ and $Ad~diag(A, B)\circ\sigma_1\circ \Pi$,
where $A, B\in GL_n(\C)$,
and $\Pi$ is defined as follows: for
$\left(
\begin{array}{cc}
a & b\\
c & d
\end{array}\right)\in sl(n,n)$, 
$\Pi \left(
\begin{array}{cc}
a & b\\
c & d
\end{array}\right)=\left(
\begin{array}{cc}
d & c\\
b & a
\end{array}\right)$  \cite{S}.
Automorphisms of the form $Ad~diag(A, B)\circ \Pi$ do not preserve the short consistent grading of $\g$, hence,
by Remark \ref{howtogetall} and Example \ref{gradedconjforfds}, every anti-linear graded conjugation of
$\g$ which is not one of those considered above, is either of the form $\varphi_1=Ad~diag(A, B)\circ \Pi\circ\sigma_1\circ\tilde{\sigma}_1$ or of the form
$\varphi_2=Ad~diag(A, B)\circ\sigma_1\circ \Pi\circ\tilde{\sigma}_1$. We have:
$\varphi_1\left(
\begin{array}{cc}
a & b\\
c & d
\end{array}\right)=
\left(
\begin{array}{cc}
A\bar{d}A^{-1} & -A\bar{c}B^{-1}\\
-B\bar{b}A^{-1} & B\bar{a}B^{-1}
\end{array}\right)$, hence
$\varphi_1^2\left(
\begin{array}{cc}
a & 0\\
0 & d
\end{array}\right)=
\left(
\begin{array}{cc}
A\bar{B}a\bar{B}^{-1}A^{-1} & 0\\
0 & B\bar{A}d\bar{A}^{-1}B^{-1}
\end{array}\right)$. It follows that if $\varphi_1$ is an anti-linear graded conjugation of $\g$, then
$A\bar{B}=\lambda I_n$ and $B\bar{A}=\sigma I_n$ for some $\lambda, \sigma\in \C$.
Besides,  
$\varphi_1^2\left(
\begin{array}{cc}
0 & b\\
0 & 0
\end{array}\right)=
\left(
\begin{array}{cc}
0 & A\bar{B}b\bar{A}^{-1}B^{-1}\\
0 & 0
\end{array}\right)$, from which it follows that $\sigma^{-1}\lambda=-1$, i.e., $\sigma=-\lambda$.
We have $A=\lambda(\bar{B})^{-1}=\lambda(\bar{\sigma})^{-1}A$, therefore $\bar{\sigma}=\lambda=-\bar{\lambda}$, i.e., 
$\lambda\in i\R$. 
By Lemma \ref{lemmasl(n,n)} $\varphi_1$ is thus equivalent either to the anti-linear graded conjugation $\tau_+$ or to
the anti-linear graded conjugation $\tau_-$, depending on $\lambda\in R^{>0}$ or $\lambda\in R^{<0}$.

Similarly, we have:
$\varphi_2\left(
\begin{array}{cc}
a & 0\\
0 & d
\end{array}\right)=
\left(
\begin{array}{cc}
A\bar{d}A^{-1} & A\bar{c}B^{-1}\\
B\bar{b}A^{-1} & B\bar{a}B^{-1}
\end{array}\right)$. Arguing as for the map $\varphi_1$ one gets statement $(b)$.

The Lie superalgebra $osp(2,2n)$ has a short consistent grading such that $\g_{\0}$ consists of matrices of the form
$\left(
\begin{array}{cc}
a & 0\\
0 & d
\end{array}\right)$, where
$a=\left(
\begin{array}{cc}
\alpha & 0\\
0 & -\alpha
\end{array}
\right)$, $\alpha\in\C$, and $d$ lies in the Lie algebra $sp(2n)$, defined by
$J_{2n}$,
 $\g_{-1}$ is the set of matrices of the form
$\left(
\begin{array}{c|c}
0 & {\begin{array}{cc}
0 & 0\\
\beta & \gamma
\end{array}}\\
\hline
{\begin{array}{cc}
\gamma^t & 0\\
-\beta^t & 0
\end{array}} & 0
\end{array}\right)$, and
$\g_{1}$ is the set of matrices of the form
$\left(
\begin{array}{c|c}
0 & {\begin{array}{cc}
\beta & \gamma\\
0 & 0
\end{array}}\\
\hline
{\begin{array}{cc}
0 & \gamma^t\\
0 & -\beta^t 
\end{array}} & 0
\end{array}\right)$, with $\beta,\gamma\in M_{1,n}$.
Every automorphism of $osp(2,2n)$ is either of the form $Ad~diag(A, B)$ for some matrices $A=diag(\alpha,\alpha^{-1})$,
$\alpha\in\C^{\times}$, $B\in Sp_{2n}(\C)$,
or of the form $Ad~diag(A, B)\circ \sigma_1$ \cite{S}. 
One can easily check that  only automorphism of the form
$Ad~diag(A, B)$ preserve the short consistent grading of $osp(2,2n)$, hence, by Remark \ref{howtogetall} and Example
\ref{gradedconjforfds}, every anti-linear graded conjugation of $osp(2,2n)$ is of the form
$\varphi=Ad~diag(A, B)\circ \tilde{\sigma}_1$ for some matrices $A=diag(\alpha,\alpha^{-1})$,
$\alpha\in\C^{\times}$, $B\in Sp_{2n}(\C)$.

Since $A$ and $a$ are diagonal matrices, we have:
$\varphi^2\left(
\begin{array}{cc}
a & 0\\
0 & d
\end{array}
\right)=
\left(
\begin{array}{cc}
a & 0\\
0 & B(\bar{B}^{-1})^td\bar{B}^tB^{-1}
\end{array}
\right)$. Therefore if $\varphi$ is an anti-linear graded conjugation of $\g$, then $B(\bar{B}^{-1})^t=\lambda I_{2n}\in Sp_{2n}$, hence 
$B=\lambda\bar{B}^t$ with $\lambda=\pm 1$, i.e., either $B$ is symplectic hermitian or it is symplectic anti-hermitian.
Besides,
$\varphi^2\left(
\begin{array}{cc}
0 & {\left(\begin{array}{cc}
0 & 0\\
\beta & \gamma
\end{array}\right)}\\
{\left(\begin{array}{cc}
\gamma^t & 0\\
-\beta^t & 0
\end{array}\right)} & 0
\end{array}
\right)=
\left(
\begin{array}{cc}
0 & -\lambda^{-1}A(\bar{A}^t)^{-1}\left({\begin{array}{cc}
0 & 0\\
\beta & \gamma
\end{array}}\right)\\
{-\lambda\left(\begin{array}{cc}
\gamma^t & 0\\
-\beta^t & 0
\end{array}\right)\bar{A}^tA^{-1}} & 0
\end{array}
\right)$. Since $A(\bar{A}^t)^{-1}=diag(\alpha\bar{\alpha}^{-1}, \alpha^{-1}\bar{\alpha})$, if $\varphi$ is an anti-linear graded conjugation of $\g$, 
either $\lambda=1$ and $\alpha=\bar{\alpha}$, i.e.\ $B$ is hermitian and $\alpha\in\R$, or $\lambda=-1$ and $\alpha=-\bar{\alpha}$, i.e.\ $B$ is anti-hermitian
and  $\alpha\in i\R$. Statement $(c)$ then follows from Lemma \ref{osp(2,2n)}. 
\hfill$\Box$

\begin{remark}\label{automorphisms}\em
If $\g$ is a simple infinite-dimensional linearly compact Lie superalgebra, then $Aut~\g$ contains a maximal reductive subgroup
which is explicitely described in \cite[Theorem 4.2]{CK2}. We shall denote this subgroup by $G$.
For any anti-linear graded conjugation $\sigma$ of $\g$, let $S$ be the subgroup of the group of linear and anti-linear
automorphisms of $\g$ generated by $\sigma$. Then
${\cal G}:=SAut~\g=Aut~\g\cup\sigma Aut~\g$ is a real algebraic group. If $\sigma$ normalizes $G$, i.e.,
$\sigma^{-1}G\sigma=G$, then $SG$ is a  subgroup of ${\cal G}$ and it is maximal reductive.
Therefore any maximal reductive subgroup
of ${\cal G}$ is conjugate into $SG$, in particular any finite order element of
${\cal G}$ is conjugate to an element of $SG$. If, moreover, $\sigma G\sigma=G$, then $SG=G\cup\sigma G$,
hence any anti-linear graded conjugation of $\g$ is conjugate to one in $\sigma G$.
\end{remark}

\begin{example}\label{K(1,2)}\em
The grading of type $(0,\dots,0|1,-1)$ of $\g=H(2k,2)$ (resp.\ $K(2k+1,2)$) is short.
Let $A=\C[[p_1,\dots,p_k,q_1,\dots,q_k]]$ (resp.\ $A=\C[[t,p_1,\dots,p_k,q_1,\dots,q_k]]$). We have:

\medskip

$\g_{-1}= \langle \xi_2\rangle\otimes A$,

$\g_0=(\langle 1,\xi_1\xi_2\rangle\otimes A)/\C 1$ (resp.\ $\langle 1,\xi_1\xi_2\rangle\otimes A$),

$\g_{1}=\langle \xi_1\rangle\otimes A$.

\medskip

\noindent
For $f\in A$, $f=\sum_{
\begin{array}{c}
(a_0, \dots, a_k)\in\Z^{k+1}_+\\
(b_1, \dots, b_k)\in\Z^k_+
\end{array}
}
\alpha_{a_0, \dots, b_k}t^{a_0}p_1^{a_1}\dots p_k^{a_k}q_1^{b_1}\dots q_k^{b_k}$, $\alpha_{a_0, \dots, b_k}\in\C$,
we denote by $\bar{f}$ the element $\bar{f}=\sum_{\begin{array}{c}
(a_0, \dots, a_k)\in\Z^k_+\\
(b_1, \dots, b_k)\in\Z^k_+
\end{array}}\overline{\alpha_{a_0, \dots, b_k}}t^{a_0}p_1^{a_1}\dots p_k^{a_k}q_1^{b_1}\dots q_k^{b_k}\in A$
(here $a_0=0$ if $A=\C[[p_1, \dots, p_k, q_1, \dots, q_k]]$). Then,
for every linear involutive change $\varphi$ of the even variables,
multiplying by $-1$ the 1-form
$\sum_{i=1}^k (p_i dq_i - q_i dp_i)$
(resp. $dt + \sum_{i=1}^k (p_i dq_i - q_i dp_i)$),
the following maps ${\Sigma}_{\varphi}^{\pm}$ are anti-linear graded conjugations of $\g$:
\begin{equation}
\begin{array}{ll}
f(p_i, q_i)\mapsto -\bar{f} ({\varphi(p_i)}, {\varphi(q_i)}) & (\mbox{resp.}~ f(t,p_i, q_i)\mapsto -\bar{f}({\varphi(t)}, {\varphi(p_i)}, {\varphi(q_i))}\\
f(p_i, q_i)\xi_1\xi_2\mapsto -\bar{f}({\varphi(p_i)}, {\varphi(q_i)})\xi_1\xi_2 & (\mbox{resp.}~ f(t,p_i, q_i)\xi_1\xi_2\mapsto -\bar{f}({\varphi(t)}, {\varphi(p_i)}, {\varphi(q_i)})\xi_1\xi_2)\\
f(p_i, q_i)\xi_1\mapsto \pm\bar{f} ({\varphi(p_i)}, {\varphi(q_i)})\xi_2 & (\mbox{resp.}~ f(t,p_i, q_i)\xi_1\mapsto \pm\bar{f}({\varphi(t)}, {\varphi(p_i)}, {\varphi(q_i)})\xi_2)\\
f(p_i, q_i)\xi_2\mapsto \mp\bar{f} ({\varphi(p_i)}, {\varphi(q_i)})\xi_1 & (\mbox{resp.}~ f(t,p_i, q_i)\xi_2\mapsto \mp\bar{f}({\varphi(t)}, {\varphi(p_i)}, {\varphi(q_i)})\xi_1).
\end{array}
\label{gradedforH}
\end{equation}
\end{example}


\begin{example}\label{S(1,2)}\em
Let $\g=S(1,2)$, $SHO(3,3)$, or $SKO(2,3;1)$. Then the algebra of outer derivations of $\g$ contains $sl_2=
\langle e,h,f\rangle$, with $e=\xi_1\xi_2\frac{\partial}{\partial x}$ and $h=\xi_1\frac{\partial}{\partial\xi_1}+\xi_2\frac{\partial}{\partial\xi_2}$ if $\g=S(1,2)$, $e=\xi_1\xi_3\frac{\partial}{\partial x_2}-\xi_2\xi_3\frac{\partial}{\partial x_1}-\xi_1\xi_2\frac{\partial}{\partial x_3}$
and $h=\sum_{i=1}^3\xi_i\frac{\partial}{\partial\xi_i}$ if
$\g=SHO(3,3)$, $e=\xi_1\xi_2\tau$ and $h=1/2(\tau-x_1\xi_1-x_2\xi_2)$ if $\g=SKO(2,3;1)$. Let us denote by $G_{out}$
the subgroup of $Aut~\g$ generated by $\exp(ad(e))$, $\exp(ad(f))$ and $\exp(ad(h))$. We recall that 
$G_{out}\subset G$, where $G$ is the subgroup of
$Aut~\g$ introduced in Remark \ref{automorphisms} \cite[Remark 2.2, Theorem 4.2]{CK2}.
We shall denote by $U_-$ the one parameter group of automorphisms $\exp(ad(tf))$,  and by
$G_{inn}$ the subgroup of $G$ consisting of inner automorphisms. 
Finally, 
$H$ will denote the subgroup of
$Aut~\g$  consisting of  invertible changes of variables multiplying the volume form
(resp.\ the even supersymplectic form) by a constant if $\g=S(1,2)$ (resp.\ $\g=SHO(3,3)$),
or the odd supercontact form by a function if $\g=SKO(2,3;1)$ (see \cite[Theorem 4.5]{CK2}). 

The gradings of type $(0|1,1)$, $(0,0,0|1,1,1)$ and $(0,0|1,1,1)$ of $\g=S(1,2)$,
$SHO(3,3)$ and $SKO(2,3;1)$, respectively, are short, and the subspaces $\g_i$'s are as follows:

\medskip

\noindent
$\g=S(1,2)$:

$\g_{-1}=\langle \frac{\partial}{\partial\xi_1}, \frac{\partial}{\partial\xi_2}\rangle\otimes\F[[x]]$

$\g_0=\{f\in \langle \frac{\partial}{\partial x}, \xi_i\frac{\partial}{\partial\xi_j} | i,j=1,2\rangle\otimes\F[[x]], div(f)=0\}$

$\g_{1}=\{f\in \langle \xi_i\frac{\partial}{\partial x}, \xi_1\xi_2\frac{\partial}{\partial\xi_i} | i=1,2\rangle\otimes\F[[x]], div(f)=0\}$.

\medskip

\noindent
$\g=SHO(3,3)$:

$\g_{-1}=\C[[x_1,x_2,x_3]]/\C 1$

$\g_0=\{f\in \langle\xi_1, \xi_2, \xi_3\rangle\otimes\C[[x_1,x_2,x_3]] | \Delta(f)=0\}$

$\g_1=\{f\in \langle\xi_i\xi_j, i,j=1,2,3\rangle\otimes\C[[x_1,x_2,x_3]] | \Delta(f)=0\}$.

\medskip

\noindent
$\g=SKO(2,3;1)$:

$\g_{-1}=\C[[x_1,x_2]]$

$\g_0=\{f\in \langle\xi_1, \xi_2, \tau\rangle\otimes\C[[x_1,x_2]] | div_1(f)=0\}$

$\g_1=\{f\in \langle\tau\xi_i, \xi_1\xi_2 ~|~ i=1,2\rangle\otimes\C[[x_1,x_2]] | div_1(f)=0\}$.

\medskip

\noindent
In all these cases the map $s=\exp(ad(e))\exp(ad(-f))\exp(ad(e))$ is a graded conjugation of $\g$: for 
$z\in\g_{-1}$, $s(z)=[e,z]$; for $z\in\g_{1}$, $s(z)=-[f,z]$, for $z\in \g_0$, $s(z)=z$.
Note that each of the above gradings can be extended to $Der ~\g=\g\rtimes {\mathfrak a}$, with ${\mathfrak a}\supset sl_2$,
so that $e$ has degree 2, $h$
has degree 0, and $f$ has degree $-2$.

We denote by $C_0$ the standard complex conjugation of $W(m,n)$ defined as follows: for 
$X=\sum_{i=1}^mP_i(x,\xi)\frac{\partial}{\partial x_i}+\sum_{j=1}^nQ_j(x, \xi)\frac{\partial}{\partial\xi_j}\in W(m,n)$,
$C_0(X)=\sum_{i=1}^m\overline{P_i(x,\xi)}\frac{\partial}{\partial x_i}+\sum_{j=1}^n\overline{Q_j(x, \xi)}\frac{\partial}{\partial\xi_j}$.
Note that $C_0$ maps $S(1,2)$ (resp.\ $SHO(3,3)$, $SKO(2,3;1)$) into itself. Moreover, $C_0\circ s=s\circ C_0$, hence, by Remark
\ref{Ccomposed}, $C_0\circ s$ is
an anti-linear graded conjugation of $S(1,2)$ (resp.\ $SHO(3,3)$, $SKO(2,3;1)$).
\end{example}

\begin{lemma}\label{SKO(2,3;beta)}
The Lie superalgebra $SKO(2,3;\beta)$
admits an anti-linear graded conjugation if and only if $|\beta|=1$, $\beta\neq\pm 1$.
\end{lemma}
{\bf Proof.} Consider the Lie superalgebra $\g=SKO(2,3;\beta)$ with its short consistent grading of type $(0,0|1,1,1)$. Then, 
by \cite[Theorem 4.2]{CK2}, there is no (linear) automorphism of $\g$ exchanging $\g_1$ and $\g_{-1}$.
Note that the restriction to $SKO(2,3;\beta)$ of the standard complex conjugation $C_0$ defined at the end of
Example \ref{S(1,2)}, defines an anti-linear Lie superalgebra isomorphism $C_{\beta}: SKO(2,3;\beta) \rightarrow SKO(2,3;\overline{\beta})$
such that  $C_{\beta}\circ C_{\bar{\beta}}=C_{\bar{\beta}}\circ C_{\beta}=id$. It follows that if $\beta\in\R$ and $\sigma$ is an anti-linear graded conjugation of $SKO(2,3;\beta)$,
then $C_\beta\circ\sigma$ is a (linear) automorphism of $SKO(2,3;\beta)$ exchanging $g_1$ and $g_{-1}$. Since such an automorphism does not exist,
we conclude that if $\beta\in\R$, then $SKO(2,3;\beta)$ has no anti-linear graded conjugations.

By \cite[Remark 4.15]{CantaK}, if $\beta\neq \pm 1$, then  $SKO(2,3;\beta)_{\bar{0}}\cong W(2,0)$ and $SKO(2,3;\beta)_{\bar{1}}\cong \Omega^0(2)^{-\frac{1}{\beta+1}}\oplus
\Omega^0(2)^{-\frac{\beta}{\beta+1}}$. Recall that a vector field $X\in W(2,0)$ acts on $\Omega^0(2)^{\lambda}$ as follows:
$X.f=X(f)+\lambda div(X)f$. If $\sigma$ is an anti-linear graded conjugation of $\g$, then $\g_{-1}$ and $\g_1$ are anti-isomorphic $\g_0$-modules,
hence $-\frac{1}{\bar{\beta}+1}=-\frac{\beta}{\beta+1}$, i.e., $|\beta|=1$. 

Now suppose $|\beta|=1$, $\beta\neq \pm 1$. Let $v$ (resp.\ $w$) be the element in $\g_{-1}$ (resp.\ $\g_1$) corresponding to $1$, and let $\sigma: \g\rightarrow \g$ be defined as
follows: 
for $X\in\g_{\bar{0}}$, $\sigma(X)=\overline{X}$; for $f\in\C[[x,y]]$, $\sigma(fv)=\bar{f}w$, $\sigma(fw)=-\bar{f}v$. Then $\sigma$ is an anti-linear graded conjugation
of $SKO(2,3;\beta)$. Indeed, one easily checks  that $\sigma([X,fv])=[\sigma(X), \bar{f}w]$.
\hfill$\Box$
\begin{proposition}\label{listofgradedconj}
The following is a complete list, up to equivalence, of anti-linear graded conjugations of 
all simple infinite-dimensional linearly
compact Lie superalgebras $\g$ with the $\Z$-gradings from Proposition \ref{fdshortconsistent}(b):
\begin{itemize}
\item[a)] $\g=H(2k,2)$: ${\Sigma}_{\varphi}^{\pm}$ are the anti-linear automorphisms of $\g$ defined by (\ref{gradedforH}).
\item[b)] $\g=K(2k+1,2)$: ${\Sigma}_{\varphi}^{\pm}$ are the anti-linear automorphisms of $\g$ defined by (\ref{gradedforH}).
\item[c)] $\g=S(1,2)$: $\tilde{\sigma}=C_0\circ s\circ \exp(ad(\alpha h))\circ \exp(ad(th'))\circ\varphi$ where $C_0$ is the standard complex conjugation, $s$ is the graded conjugation of $\g$ introduced in Example \ref{S(1,2)}, $t\in i\R$, $h'=2x\frac{\partial}{\partial x}+\xi_1\frac{\partial}{\partial \xi_1}+\xi_2\frac{\partial}{\partial \xi_2}$, $\varphi$ is an element of the $SL_2$-subgroup of $G_{inn}$ generated by $exp(ad(\xi_1\frac{\partial}{\partial \xi_1}-\xi_2\frac{\partial}{\partial \xi_2}))$,
$exp(ad(\xi_1\frac{\partial}{\partial \xi_2}))$ and $exp(ad(\xi_2\frac{\partial}{\partial \xi_1}))$   such that
$\overline{\varphi}\varphi=I_2$ (resp.\ $\overline{\varphi}\varphi=-I_2$), and $\exp(\alpha)=\pm 1$ (resp.\ $\exp(\alpha)=\pm i$). 
\item[d)] $\g=SHO(3,3)$: $\tilde{\sigma}=C_0\circ s\circ\exp(ad(\alpha h))\circ\varphi$ with $\varphi$ such that $\overline{\varphi}\varphi=1$
(resp.\ $\overline{\varphi}\varphi=-1$), lying in the subgroup of $Aut \g$ generated by  $G_{inn}$ and 
$\exp(ad(\Phi))$, $\Phi=\sum_{i=1}^n(-x_i\frac{\partial}{\partial x_i}+\xi_i\frac{\partial}{\partial \xi_i})$,
$\exp(\alpha)=\pm 1$ (resp.\ $\exp(\alpha)=\pm i$).
\item[e)] $\g=SKO(2,3;1)$: 
$\tilde{\sigma}=C_0\circ s\circ\exp(ad(\alpha h))\circ\varphi$, where $C_0$ is the standard complex conjugation, $s$ is the graded conjugation of $\g$ introduced in Example \ref{S(1,2)}, $\exp(\alpha)=\pm 1$  and $\varphi$ is an element in $G_{inn}$ such that $\overline{\varphi}\varphi=I_2$.
\item[f)] $\g=SKO(2,3;\beta)$,  $|\beta|=1$, $\beta\neq \pm 1$: $\tilde{\sigma}=\sigma\circ\exp(ad(\alpha h))\circ\varphi$, where $\sigma$ is defined as in Lemma \ref{SKO(2,3;beta)}, $h=\frac{\tau-x_1\xi_1-x_2\xi_2}{2}$, $\varphi$ is an element in $G_{inn}$ such that
$\overline{\varphi}\varphi=I_2$ and $\exp(\alpha)=\pm 1$.
\end{itemize}
\end{proposition}
{\bf Proof.} 
Let $\g=H(2k,2)$ with the grading of type $(0,\dots,0|1,-1)$ (see Example \ref{K(1,2)}), let
$\Sigma_{\varphi}$ be the graded conjugation of $\g$ defined in Example \ref{K(1,2)}, and let $\sigma$
be an anti-linear graded conjugation of $\g$. By \cite[Theorem 4.2]{CK2}, the group
$G=\C^{\times}(Sp_{2k}\times O_2)$ consists of linear changes of variables preserving the symplectic form
up to multiplication by a non-zero scalar \cite[Theorem 4.2]{CK2}, hence the map $\Sigma_{\varphi}$ satisfies both conditions
$\Sigma_{\varphi}^{-1}G\Sigma_{\varphi}=G$ and $\Sigma_{\varphi}G\Sigma_{\varphi}=G$.
Then, by Remark \ref{howtogetall} and Remark \ref{automorphisms}, we may assume, up to conjugation, that $\sigma=F\circ\Sigma_{\varphi}$
where  $F$ is a grading preserving automorphism of $\g$ such that 
$(F\circ\Sigma_{\varphi})^2(x_j)=(-1)^j x_j$ for $x_j\in\g_j$, and
$F\in G$.  

Since $F$ preserves the grading and lies in $G$,
we have: $F(\xi_1)=a\xi_1$ and $F(\xi_2)=b\xi_2$, for some $a,b\in\C^\times$. 
Besides, the condition 
$(F\circ\Sigma_{\varphi})^2(\xi_i)=-\xi_i$ implies $b=(\bar{a})^{-1}$.

It follows that, for $f\in\C[[p_i,q_i]]$,
$F(f)=F([\xi_1, f\xi_2])=[a\xi_1, F(f\xi_2)]=[a\xi_1, \tilde{f}\xi_2]=a\tilde{f}$ for some
$\tilde{f}\in\C[[p_i,q_i]]$, i.e., $F(f\xi_2)=a^{-1}F(f)\xi_2$. Likewise,
$F(f\xi_1)=\bar{a}F(f)\xi_1$  and $F(f\xi_1\xi_2)=F(f)\xi_1\xi_2$. Note that 
$F([\xi_1\xi_2, \xi_1]=F(\xi_1)=a\xi_1$ and
$F([\xi_1\xi_2, \xi_1]=[F(\xi_1\xi_2), a\xi_1]=a^2\bar{a}^{-1}\xi_1$, hence $a\in \R$.

It follows that $F\circ\Sigma_{\varphi}$ is defined as follows:
\begin{equation}
\begin{array}{l}
f(p_i, q_i)\mapsto -\bar{f} (F\circ\varphi(p_i), F\circ\varphi(q_i))\\
f(p_i, q_i)\xi_1\xi_2\mapsto -\bar{f} (F\circ\varphi(p_i), F\circ\varphi(q_i))\xi_1\xi_2\\
f(p_i, q_i)\xi_1\mapsto a\bar{f} (F\circ\varphi(p_i), F\circ\varphi(q_i))\xi_2\\
f(p_i, q_i)\xi_2\mapsto -a^{-1}\bar{f} (F\circ\varphi(p_i), F\circ\varphi(q_i))\xi_1
\end{array}
\end{equation}
for some linear change of even variables $\varphi$, such that $(F\circ \varphi)^2=1$.
Besides,
since $F(\xi_1)=a\xi_1$
and $F(\xi_2)=a^{-1}\xi_2$,
$F$ preserves the odd part $d\xi_1d\xi_2$ of the symplectic form,
hence it preserves  the symplectic form.
Finally, by Remark \ref{scalars} we may assume $a=\pm 1$. 
This concludes the proof of a). Similar arguments prove b).

Let $\g=S(1,2)$. We already noticed in Example \ref{S(1,2)} that $\tilde{s}=C_0\circ s$ is an anti-linear graded conjugation of $\g$.
Here $G_{inn}$ is generated by $exp(ad(h'))$, $exp(ad(\xi_1\frac{\partial}{\partial \xi_1}-\xi_2\frac{\partial}{\partial \xi_2}))$,
$exp(ad(\xi_1\frac{\partial}{\partial \xi_2}))$ and $exp(ad(\xi_2\frac{\partial}{\partial \xi_1}))$ and $G=\C^{\times}SO_4$ is generated by
$G_{inn}$ and the automorphisms $\exp(ad(e))$, $\exp(ad(f))$ and $\exp(ad(h))$ (see Example \ref{S(1,2)} and \cite[Theorem 4.2]{CK2}). Therefore $\tilde{s}$ satisfies both relations
$\tilde{s}^{-1}G\tilde{s}=G$ and $\tilde{s}G\tilde{s}=G$, i.e, it satisfies the hypotheses of Remark \ref{automorphisms}, hence
every anti-linear graded conjugation of $\g$ is conjugate to one in $\tilde{s}G$. By Remark \ref{howtogetall} we therefore 
aim to classify all automorphisms $\psi\in G$  preserving the grading of $\g$
of type $(0|1,1)$, such that $(\tilde{s}\circ\psi)^2(x_j)=(-1)^jx_j$ for $x_j\in\g_j$.

By \cite[Remark 4.6]{CK2}, if $\psi$ is an automorphism of $\g$ lying in $G$, then either $\psi\in U_-H\cap G$ or $\psi\in U_-sH\cap G$.
Note that $U_{-}H\cap G=U_-(H\cap G)$ and $U_{-}sH\cap G=U_-s(H\cap G)$, since $U_{-}\subset G$ and $s\in G$ \cite[Theorem 4.2]{CK2}. 
Here $H\cap G$ is the subgroup of $Aut~\g$ generated by $\exp(ad(e))$, $\exp(ad(h))$ and $G_{inn}$. Note that
$G_{inn}\subset \exp(ad(\g_0))$.

Let $\psi\in U_-(H\cap G)$.  
Then $\psi=\exp(ad(tf))\psi_0$ for some $t\in\C$ and some $\psi_0\in H\cap G$. For $x\in\g_1$, we have:
$\psi(x)=\exp(ad(tf))(\psi_0(x))=\psi_0(x)+t[f,\psi_0(x)]$, since $\psi_0(x)\in \g_1$. Since $\psi$ preserves the grading, $t=0$, i.e., 
$\psi=\psi_0\in H\cap G$. We may hence write $\psi=\varphi_0\varphi_1$ for some
$\varphi_1\in G_{inn}$ and some $\varphi_0$ lying in the subgroup generated by
$\exp(ad(h))$ and $\exp(ad(e))$, and assume
$\varphi_0=\exp(ad(\beta e))\exp(ad(\alpha h))$ for some $\alpha, \beta\in\C$.
For $z\in\g_{-1}$ we have: $\psi(z)=\varphi_0(\varphi_1(z))=exp(-\alpha)\exp(ad(\beta e))(\varphi_1(z))=
\exp(-\alpha)(\varphi_1(z)+\beta[e,\varphi_1(z)])$, since $\varphi_1(z)\in \g_{-1}$. Since $\psi$ preserves
the grading, we have $\beta=0$, i.e., $\psi=\exp(ad(\alpha h))\varphi_1$, hence
$$\tilde{\sigma}=\tilde{s}\circ\psi=C_0\circ s\circ\exp(ad(\alpha h))\circ\varphi_1.$$
We have:
$$\varphi_1\circ C_0=C_0\circ\overline{\varphi_1},$$
where, for $\varphi=\prod_{i=1}^k\exp(ad(z_i))$, $\overline{\varphi_1}=\prod_{i=1}^k\exp(ad(C_0(z_i)))$, and, for $\alpha\in\C$,
$$\exp(ad(\alpha h)) \circ C_0=C_0\circ \exp(ad(\bar{\alpha}h)),\,\,\,\,\,\exp(ad(\alpha h)) \circ s=s\circ \exp(ad(-{\alpha}h)).$$
It follows that $(\tilde{\sigma})^2=(\tilde{s})^2\circ\exp(ad(\alpha-\bar{\alpha})h)\circ\overline{\varphi_1}\varphi_1$,
hence 
\begin{equation}
\exp(ad(\alpha-\bar{\alpha})h)\circ\overline{\varphi_1}\varphi_1=id.
\label{C}
\end{equation}
In particular, 
$\exp(ad(\alpha-\bar{\alpha})h)\circ\overline{\varphi_1}\varphi_1(\frac{\partial}{\partial x})=\frac{\partial}{\partial x}$.
Since $\varphi_1$ lies in $G_{inn}$, we may assume that $\varphi_1=\exp(ad(th'))\varphi_2=\varphi_2\exp(ad(th'))$ for some
$\varphi_2$ in the $SL_2$-subgroup of $G_{inn}$ 
generated by $exp(ad(\xi_1\frac{\partial}{\partial \xi_1}-\xi_2\frac{\partial}{\partial \xi_2}))$,
$exp(ad(\xi_1\frac{\partial}{\partial \xi_2}))$ and $exp(ad(\xi_2\frac{\partial}{\partial \xi_1}))$, therefore
$\overline{\varphi_1}\varphi_1=\exp(ad(t+\bar{t})h')\overline{\varphi_2}\varphi_2$. Since $\overline{\varphi_2}\varphi_2
(\frac{\partial}{\partial x})=\frac{\partial}{\partial x}$, and $\exp(ad(\alpha-\bar{\alpha})h)(\frac{\partial}{\partial x})=\frac{\partial}{\partial x}$ we have:
$\exp(ad(t+\bar{t})h')(\frac{\partial}{\partial x})=\exp(-2(t+\bar{t}))(\frac{\partial}{\partial x})$, hence $t\in i\R$.

If we restrict condition (\ref{C}) to the $sl_2$-subalgebra of $\g_0$ generated by $\xi_1\frac{\partial}{\partial \xi_2}$ and
$\xi_2\frac{\partial}{\partial \xi_1}$, we find that $\overline{\varphi_2}\varphi_2$ lies in $\langle I_2\rangle\cap SL_2$, i.e.,
$\overline{\varphi_2}\varphi_2=\pm I_2$.
Then, if we restrict condition (\ref{C}) to the subspace $S=\langle \frac{\partial}{\partial\xi_1}, \frac{\partial}{\partial\xi_2}\rangle\subset \g_{-1}$,
we get:
$$\exp(ad(\alpha-\bar{\alpha})h)\exp(ad(t+\bar{t})h') \overline{\varphi_2}\varphi_2|_S=id,$$
i.e., either $\overline{\varphi_2}\varphi_2=I_2$ and $\exp(-(\alpha-\bar{\alpha}))=1$, or $\overline{\varphi_2}\varphi_2=-I_2$
and $\exp(-(\alpha-\bar{\alpha}))=-1$. Using Remark \ref{scalars}, we get statement $(c)$.

Finally note that, since $\psi$ preserves the grading, it cannot lie in $U_{-}s(H\cap G)$. Indeed, suppose that $\psi\in U_{-}s(H\cap G)$.
Then we may assume that $\psi=\exp(ad(tf))\circ s\circ\exp(ad(\beta e))\exp(ad(\alpha h))\varphi_1$, for some $\alpha, \beta, t\in\C$ and some $\varphi_1\in G_{inn}$.
Suppose that $x\in\g_1$. Then $\varphi_1(x)\in \g_1$, since $G_{inn}\subset \exp(ad(\g_0))$, hence $\exp(ad(\beta e))(\varphi_1(x))=\varphi_1(x)$. It follows that
$\psi(x)=\exp(\alpha)\exp(ad(tf))(s(\varphi_1(x)))$ $=\exp(\alpha)s(\varphi_1(x))$, since $s|_{\g_1}:\g_1 \rightarrow \g_{-1}$. Therefore
$\psi|_{\g_1}:\g_1 \rightarrow \g_{-1}$, a contradiction.

Let $\g=SHO(3,3)$ with the grading of type $(0,0,0|1,1,1)$. 
We already noticed in Example \ref{S(1,2)} that $\tilde{s}=C_0\circ s$ is an anti-linear graded conjugation of $\g$.
Here $G_{inn}$ is generated by automorphisms $exp(ad(x_i\xi_j))$ with $i\neq j$, $\exp(ad(e))$, $\exp(ad(f))$, $\exp(ad(h))$ and $\exp(ad(\Phi))$ (see Example \ref{S(1,2)} and \cite[Theorem 4.2]{CK2}). Therefore $\tilde{s}$ satisfies both relations
$\tilde{s}^{-1}G\tilde{s}=G$ and $\tilde{s}G\tilde{s}=G$, i.e, it satisfies the hypotheses of Remark \ref{automorphisms}, hence
every anti-linear graded conjugation of $\g$ is conjugate to one in $\tilde{s}G$. By Remark \ref{howtogetall} we therefore 
need to classify all automorphisms $\psi\in G$  preserving the grading of $\g$
of type $(0,0,0|1,1,1)$, such that $(\tilde{s}\circ\psi)^2(x_j)=(-1)^jx_j$ for $x_j\in\g_j$.

By \cite[Remark 4.6]{CK2}, if $\psi$ is an automorphism of $\g$ lying in $G$, then either $\psi\in U_-H\cap G$ or $\psi\in U_-sH\cap G$.
Like in the case of $S(1,2)$, we have: $U_{-}H\cap G=U_-(H\cap G)$ and $U_{-}sH\cap G=U_-s(H\cap G)$, since $U_{-}\subset G$ and $s\in G$ \cite[Theorem 4.2]{CK2}. 
Here $H\cap G$ is the subgroup of $Aut~\g$ generated by $\exp(ad(e))$, $\exp(ad(h))$, $\exp(ad(\Phi))$ and $G_{inn}$. Note that
$G_{inn}\subset \exp(ad(\g_0))$.

Let $\psi\in U_-(H\cap G)$.  
Then $\psi=\exp(ad(tf))\psi_0$ for some $t\in\C$ and some $\psi_0\in H\cap G$. For $x\in\g_1$, we have:
$\psi(x)=\exp(ad(tf))(\psi_0(x))=\psi_0(x)+t[f,\psi_0(x)]$, since $\psi_0(x)\in \g_1$. Since $\psi$ preserves the grading, $t=0$, i.e., 
$\psi=\psi_0\in H\cap G$. We may hence write $\psi=\exp(ad(\gamma e))\varphi$ for some $\gamma\in \C$ and some
$\varphi$ in the subgroup of $Aut \g$ generated by $\exp(ad(h))$, $\exp(ad(\Phi))$ and $G_{inn}$.
For $z\in\g_{-1}$ we have: $\psi(z)=\varphi(z)+\gamma[e,\varphi(z)]$, since $\varphi(z)\in \g_{-1}$. Since $\psi$ preserves
the grading, we have $\gamma=0$, i.e., $\psi=\exp(ad(\beta\Phi))\exp(ad(\alpha h))\varphi_0$ for some $\alpha, \beta\in\C$ and $\varphi_0\in G_{inn}$, hence
every anti-linear graded conjugation $\tilde{\sigma}$ of $\g$ is of the form:
$$\tilde{\sigma}=\tilde{s}\circ\psi=C_0\circ s\circ\exp(ad(\beta\Phi))\circ\exp(ad(\alpha h))\circ\varphi_0.$$
We have:
$$\varphi_0\circ C_0=C_0\circ\overline{\varphi_0},$$
where, for $\varphi=\prod_{i=1}^k\exp(ad(z_i))$, $\overline{\varphi}=\prod_{i=1}^k\exp(ad(C_0(z_i)))$. Besides, for $\alpha, \beta\in\C$,
$$\exp(ad(\alpha h)) \circ C_0=C_0\circ \exp(ad(\bar{\alpha}h)),\,\,\,\,\,\exp(ad(\alpha h)) \circ s=s\circ \exp(ad(-{\alpha}h));$$
$$\exp(ad(\beta \Phi)) \circ C_0=C_0\circ \exp(ad(\bar{\beta}h)),\,\exp(ad(\beta\Phi)) \circ s|_{\g_k}=exp(-3k\beta)s\circ \exp(ad({\beta}\Phi))|_{\g_k}\,(k=-1,0,1).$$
It follows that $(\tilde{\sigma})^2|_{\g_k}=exp(-3k\bar{\beta})(\tilde{s})^2\circ\exp(ad(\beta+\bar{\beta})\Phi)\circ\exp(ad(\alpha-\bar{\alpha})h)\circ\overline{\varphi_0}\varphi_0|_{\g_k}$,
hence, since $\varphi_0$, $\overline{\varphi_0}$, $\exp(ad(h))$ and $\exp(ad(\Phi))$ preserve the subspace $\g_k$ for every $k$, 
\begin{equation}
exp(-3k\bar{\beta})\exp(ad(\beta+\bar{\beta})\Phi)\circ\exp(ad(\alpha-\bar{\alpha})h)\circ\overline{\varphi_0}\varphi_0|_{\g_k}=id\,\, (k=-1,0,1).
\label{D}
\end{equation}

Let us restrict condition (\ref{D}) to the subalgebra $S$ of $\g_0$ generated by the elements $x_i\xi_j$ with $1\leq i\neq j\leq 2$: $S\cong sl_3$. Then
$G_{inn}$ act on $S$ via the Adjoint action and $\exp(ad(h))|_S=id_S=\exp(ad(\Phi))|_S$. Hence (\ref{D}) becomes simply  
 $\overline{\varphi_0}\varphi_0|_S=id$, i.e.,
$\overline{\varphi_0}\varphi_0=\lambda I_3$ for $\lambda\in\C$ such that $\lambda^3=1$.

Now let us restrict condition (\ref{D}) to the subspace $W$ of $\g_0$ spanned by the elements $\xi_i$, for $i=1,2,3$. Then $G_{inn}$ acts on $W$ via the standard
action and $\exp(ad(h))|_W=id_W$, $\exp(ad(\beta\Phi))|_W=\exp(\beta)id_W$. Hence (\ref{D}) becomes: $\lambda \exp(\beta+\bar{\beta})=1$. Since $\lambda^3=1$, it follows
that $\lambda=1$ and $\beta\in i\R$. Therefore (\ref{D}) reduces to the following relation:
\begin{equation}
exp(-3k\bar{\beta})\exp(ad(\alpha-\bar{\alpha})h)\circ\overline{\varphi_0}\varphi_0|_{\g_k}=id\,\, (k=-1,0,1).
\label{Dreduced}
\end{equation}
By Remark \ref{scalars}, we can assume $\alpha\in\{i\pi/2, i\pi, i3\pi/2, 2i\pi\}$, i.e., $\exp(\bar{\alpha}-\alpha)=\pm 1$. It follows that if we restrict
condition (\ref{Dreduced}) to the subspace $V=\langle x_1, x_2, x_3\rangle\subset\g_{-1}$, we find that either $\exp(3\bar{\beta})=1$ if
$\exp(\bar{\alpha}-\alpha)=1$, or  $\exp(3\bar{\beta})=-1$ if
$\exp(\bar{\alpha}-\alpha)=-1$. Hence $\varphi:=\exp(ad(\beta\Phi))\circ\varphi_0$ satisfies statement $d)$.

Finally note that, since $\psi$ preserves the grading, it cannot lie in $U_{-}s(H\cap G)$. Indeed, suppose that $\psi\in U_{-}s(H\cap G)$.
Then we may assume that $\psi=\exp(ad(tf))\circ s\circ\varphi$, for some $t\in\C$ and some $\varphi\in H\cap G$.
Suppose that $x\in\g_1$. Then $\varphi(x)\in \g_1$, hence $s(\varphi(x))\in \g_{-1}$. It follows that
$\psi(x)=\exp(ad(tf))(s(\varphi(x)))\in\g_{-1}$. Therefore
$\psi|_{\g_1}:\g_1 \rightarrow \g_{-1}$, a contradiction.  

The argument for $\g=SKO(2,3;1)$ and $\g=SKO(2,3;\beta)$ with $\beta\neq 1$ is the same as for $\g=S(1,2)$. 
\hfill$\Box$

\end{document}